\pdfoutput=1
\documentclass[10pt]{article}
\usepackage[a4paper, margin=1in]{geometry}
\usepackage{changepage}
\usepackage{amsmath, amssymb, amsthm}
\usepackage[utf8]{inputenc}
\usepackage[english]{babel}
\usepackage[numbers]{natbib}
\usepackage{url}
\usepackage[usenames]{color}
\usepackage[colorlinks=true]{hyperref}
\usepackage{cleveref}
\usepackage{tabularx}
\usepackage[title]{appendix}
\usepackage{xcolor}
\usepackage{graphicx}
\usepackage{placeins}
\usepackage{cancel}
\theoremstyle{plain}

\newtheorem{lemma}{Lemma}[section]

\newtheorem{theorem}{Theorem}[section]
\newtheorem*{theorem*}{Theorem}
\newtheorem{remark}{Remark}[section]

\theoremstyle{definition}
\newtheorem{question}{Question}
\crefname{conjecture}{Conjecture}{Conjectures}
\crefname{theorem}{Theorem}{Theorems}
\crefname{theorem*}{Theorem}{Theorems}
\crefname{corollary}{Corollary}{Corollaries}
\crefname{lemma}{Lemma}{Lemmas}
\crefname{proposition}{Proposition}{Propositions}
\crefname{remark}{Remark}{Remarks}
\crefname{note}{Note}{Notes}
\crefname{definition}{Definition}{Definitions}
\crefname{notation}{Notation}{Notations}
\crefname{example}{Example}{Examples}
\crefname{question}{Question}{Questions}
\crefname{section}{\S}{Sections}
\crefname{equation}{Equation}{Equations}
\crefformat{equation}{(eq.~#2#1#3)}

\newcommand{\floor}[1]{\left\lfloor #1 \right\rfloor}

\newcommand{\numerator}{\operatorname{Numerator}}

\newcommand{\HW}{\operatorname{HW}}

\newcommand{\Numer}{L}
\newcommand{\Denom}{D}
\newcommand{\ATerm}{\mathcal{A}}
\newcommand{\CTerm}{\mathcal{C}}
\newcommand{\NTerm}{\mathcal{N}}

\newcommand{\Z}{\mathbb{Z}}
\newcommand{\Kalmar}{\mathcal{E}^3}
\DeclareMathOperator{\DTIME}{DTIME}
\DeclareMathOperator{\ELEMENTARY}{ELEMENTARY}

\newcommand{\seqnum}[1]{\href{https://oeis.org/#1}{\rm \underline{#1}}}

\title{On arithmetic terms expressing the prime-counting function and the n-th prime}
\author{Mihai Prunescu \footnote{Research Center for Logic, Optimization and Security (LOS), Faculty of Mathematics and Computer Science, University of Bucharest, Academiei 14, Bucharest (RO-010014), Romania; and Simion Stoilow Institute of Mathematics of the Romanian Academy, Research unit 5, P. O. Box 1-764, Bucharest (RO-014700), Romania. E-mail: {\tt mihai.prunescu@imar.ro}, {\tt mihai.prunescu@gmail.com}.},
Joseph M. Shunia \footnote{Wraithwatch, Austin, Texas, United States. E-mail: {\tt jshunia@gmail.com}.}}
\date{December 2024 \\ \footnotesize{Revised: August 2025}}

\begin{document}

\maketitle

\begin{abstract} \noindent
We present the first fixed-length elementary closed-form expressions for the prime-counting function, $\pi(n)$, and the $n$-th prime number, $p(n)$. These expressions are arithmetic terms, requiring only a finite and fixed number of elementary arithmetic operations from the set: addition, subtraction, multiplication, integer division, and exponentiation.
\\[2mm]
Mazzanti proved that every Kalmar function can be represented as an arithmetic term. We develop an arithmetic term representing the prime omega function, $\omega(n)$, which counts the number of distinct prime divisors of a positive integer $n$. From this term, we find immediately an arithmetic term for the prime-counting function, $\pi(n)$. Combining these results with a new arithmetic term for binomial coefficients and novel prime-related exponential Diophantine equations, we manage to develop an arithmetic term for the $n$-th prime number, $p(n)$, thereby providing a constructive solution to the fundamental question: Is there an order to the primes?
\\[2mm]
\end{abstract}

\section{Introduction}
The \textbf{prime numbers} are the most fundamental elements in arithmetic, as every natural number greater than $1$ can be expressed as the product of one or more primes. Primes have been studied for millennia, since their treatment arises naturally from the study of counting and multiplication.

The \textbf{$n$-th prime number} is represented as $p(n)$ and the sequence of prime numbers begins as:
\begin{align*}
2, 3, 5, 7, 11, 13, 17, 19, 23, 29, 31, 37, 41, 43, 47, 53, 59, 61, 67, 71, 73, 79, 83, 89, 97, 101, \ldots
\quad \text{(see \seqnum{A000040})}
\end{align*}
The primes display an erratic growth and upon initial inspection, appear to be distributed randomly. Yet, a deeper analytical approach suggests a hidden order. This enigmatic order is evident in the prime number theorem, which states $\pi(n) \sim \frac{n}{\log(n)}$, where $\pi(n)$ represents \textbf{the prime-counting function}, returning the number of primes less than or equal to $n$.

The primes are essential in mathematics and defined with remarkable simplicity, yet exhibit a chaotic and seemingly unpredictable growth. This duality has inspired extensive study, shaping the development of number theory itself, where primes remain a central focus \cite{dusautoy2012musicoftheprimes}. Among the most profound questions in this field is the truth or falsehood of the Riemann hypothesis, which conjectures that all non-trivial zeros of the Riemann zeta function $\zeta(z) = \sum_{n=1}^\infty \frac{1}{n^z}$ lie on the critical line $\Re(z) = \frac{1}{2}$. Its resolution is widely regarded as one of the most important unsolved problems in mathematics, promising to provide insight into the intricate structure underlying the distribution of prime numbers \cite{devlin2001mathematics}.

\subsection{The order of the primes}
\begin{quote}
\textit{``Mathematicians have tried in vain to this day to discover some order in the sequence of prime numbers, and we have reason to believe that it is a mystery into which the human mind will never penetrate.''} -- Leonhard Euler, 1770 \cite{simmons2007calculusgems}
\end{quote}

A most fundamental question is if there exists an \textbf{order} to the prime numbers: A deterministic pattern or structure within the natural numbers that dictates the sequence of primes. Formally, one might seek a function $f : \mathbb{N} \to \mathbb{N}$ defined using only a finite and fixed number of elementary arithmetic operations from the set $\{+, -, \cdot, /, x^y\}$, which for a given $n$, returns the $n$-th prime number. Such a function, if it exists, would encapsulate any inherent order or pattern in the primes. Yet constructing such a function has proven to be an extraordinary challenge.

To illustrate why, consider the classical \textbf{Sieve of Eratosthenes}, an algorithm for generating prime numbers. Beginning with the natural numbers greater than $1$, the sieve starts at the first prime $p(1)=2$ and crosses out all its multiples. It then moves to next unmarked number $p(2)=3$, crossing out its multiples. The process of advancing and crossing-out continues indefinitely, leaving only the prime numbers unmarked. This sieve is inherently sequential and iterative, as the determination of the $n$-th prime depends not only on the previous prime, but on \textit{all} primes that precede it. The sequence of primes is infinite. So then, how could one possibly capture this behavior in a finite and fixed-length formula?

Indeed, the problem of discerning a precise order for the prime numbers has remained entirely out of reach. That is, until relatively recent advancements in mathematical logic and computability unveiled a promising approach.

\subsection{Kalmar functions}

The \textbf{Kalmar functions}, also called \textbf{elementary functions}, are the class $\Kalmar$ in the Grzegorczyk hierarchy \cite{grzegorczyk1953someclasses}, consisting of all primitive recursive functions $f : \mathbb{N}^k \to \mathbb{N}$ that can be computed deterministically in iterated exponential time \cite{marchenkov2007superposition, oitavem1997}. Specifically, there exists a constant $d \in \mathbb{N}$ such that for any input $\vec{n} = (n_1, n_2, \ldots, n_k) \in \mathbb{N}^k$, $f(\vec{n})$ can be computed in time $O(\exp_d(\|\vec{n}\|_{\infty}))$, where $\exp_d(\cdot)$ denotes the \textbf{$d$-fold exponential function} and $\|\vec{n}\|_{\infty} := \max(n_1,n_2,\ldots,n_k)$ represents the \textbf{infinity norm} or largest element of $\vec{n}$. The corresponding decision class for $\Kalmar$ is:
\begin{align*}
\ELEMENTARY = \bigcup_{d \in \mathbb{N}} \DTIME(\exp_d(\|\vec{n}\|_{\infty})) .
\end{align*}

\textbf{Arithmetic terms} are defined in \cite{mazzanti2002plainbases, marchenkov2007superposition} as functions $f : \mathbb{N}^k \to \mathbb{N}$ that can be computed deterministically in iterated exponential time and expressed in the language:
\begin{align*}
L = \{ +, \dot{-}, \cdot, /, x^y \} ,
\end{align*}
where the \textbf{monus} operator ($\dot{-}$) denotes bounded subtraction, defined as: $a \dot{-} b = \max(a - b, 0)$ (see \citep[pg.~141]{vereschchaginshen2002computable}). The use of bounded subtraction ensures that outputs remain in $\mathbb{N}$, though standard subtraction ($-$) may be used equivalently, provided that the output remains a natural number. The ($/$) operator denotes integer division and the $\bmod$ operation is implicitly included in the set, since $a \bmod b = a - b \floor{a/b}$. It is important to emphasize that all arithmetic terms are \textbf{elementary closed-form expressions} and, more restrictively, are of \textbf{fixed-length}, meaning they do not permit variable-length summations and products like those typically allowed in elementary closed-form expressions. Precisely, arithmetic terms are \textbf{fixed-length elementary closed-form expressions}. There is no unanimous consensus of what constitutes a \textbf{closed-form}, however see \cite{chow1999closedformnumber} for an attempt at a rigorous definition.

\subsection{On constructing arithmetic terms for primes}

Mazzanti proved in \cite{mazzanti2002plainbases} that every Kalmar function can be represented by arithmetic terms (see also Marchenkov \cite{marchenkov2007superposition}). This is a rather surprising result, since the function $p(n)$, computing the $n$-th prime number, is primitive recursive and bounded above by $O(n^2)$ \cite{jones1975formula}. Hence, $p(n)$ is a Kalmar function. The subtle implication being: \textit{There exists an arithmetic term to compute the $n$-th prime number.}

However, Mazzanti's result on this matter is purely theoretical and non-constructive, leaving open the possibility that any such arithmetic term for $p(n)$ might be so large that it could not be practically realized. Mazzanti's approach, the \textbf{hypercube method} (described in \cref{section:hypercube}), makes clever usage of elementary arithmetic to count the number of solutions to Diophantine equations. Hilbert's 10th problem, which asked for a general algorithm that can determine if an arbitrary Diophantine equation has solutions in $\mathbb{N}$, was shown to be unsolvable by Matiyasevich (for the details, see \cite{matiyasevich1993hilbert}). Thus, while the hypercube method offers a potential approach, there can be no general algorithm nor procedure for constructing arithmetic terms by its application.

As will be demonstrated, this task of constructing an arithmetic term for $p(n)$ is far from straightforward; in addition to inherit challenges with the hypercube method we described above, extremely large computations and equations are often required, making proofs necessarily complicated and potentially out of reach. In the case of $p(n)$, the mathematical formulas and symbolic computations involved can explode in size, rapidly become unwieldy.

The difficulty in constructing an arithmetic term for $p(n)$ is such that it necessitated two complete rewrites of this paper. The initial version was based on \textbf{Wilson's theorem}, which states that $n$ is prime if and only if $(n-1)!^2 \equiv 1 \pmod{n}$. In this first version, the computations for $p(n)$ became so immense that it was exceedingly complicated to derive certain bounds required to obtain the final arithmetic term and complete the proof. Our subsequent discovery of an arithmetic term for the \textbf{prime omega function} $\omega(n)$, which counts the number of distinct prime divisors of $n$ (see \cref{section:omega}), provided massive simplifications and an entirely new approach to $p(n)$. While rewriting the paper, we were able to prove a new arithmetic term for binomial coefficients $\binom{a}{b}$ (see \cref{section:binomialcoefficients}), yielding another significant simplification to the arithmetic term for $p(n)$ and its proofs, culminating in this third and final version. In the end, we found a way to apply the hypercube method to obtain an arithmetic term for $p(n)$ that is massive and computationally impractical, but can actually be written down and proved.

\subsection{Formulas for primes} \label{subsection:primeformulas}
A notable formula for computing the $n$-th prime was introduced by Willans in \cite{willans1964formula}:
\begin{align} \label{eq:willans}
p(n) = 1 + \sum_{i=1}^{2^n} \floor{\left(\frac{n}{\sum_{j=1}^i \floor{\left(\cos \frac{(j-1)! + 1}{j} \pi\right)^2} }\right)^{1/n}} .
\end{align}
Although this formula provides an explicit expression for determining primes, it is not an arithmetic term due to its reliance on summations of variable length. It operates by performing an obfuscated primality test on each number in the sum, based on Wilson's theorem. The key component, $\cos(\frac{(j-1)! + 1}{j} \pi)^2$, encodes the $(n-1)!^2 \bmod \; n$ operation using trigonometric functions, evaluating to $1$ if $j$ is prime and $0$ otherwise. This result is then used to locate $p(n)$. There are some well-known variants of Willans formula, such as that of Jones \cite{jones1975formula}. It is worth mentioning that Jones' formula is technically an elementary closed-form, though like Willans' formula, it contains a variable length summation and is therefore also not an arithmetic term. While clever, such formulas rely previous knowledge of the primality of $p(n)$; hence, they are entirely self-referential, and deduce $p(n)$ only by circular means.

An arithmetic term for the factorial function, $n!$, was first discovered by Robinson in \cite{robinson1952arithmetic}. Unaware of Robinson's earlier result, Prunescu and Sauras-Altuzarra published another version in \cite{prunescu2024factorial} where one also finds as an application, an arithmetic term whose image is the set of primes, but with repetitions:
\begin{align} \label{eq:prunescusaurasprimefunction}
f(n) = 2 + \left( (2 \cdot n!) \bmod (n+1) \right)
= \begin{cases}
n + 1 \quad \textup{ if } n + 1 \textup{ is prime} \\
2 \quad \textup{ otherwise}
\end{cases} ,
\end{align}
where the arithmetic term for $n!$ is given as
\begin{align}
n! = \floor
{
    \frac
    {
        2^{n(n+1)(n+2)}
    }
    {
        \floor
        {
            \left(
                2^{2^{(n+1)(n+2)}-n} + 2^{-n}
            \right)^{2^{(n+1)(n+2)}}
        }
        \bmod
        2^{2^{(n+1)(n+2)}}
    }
} ,
\end{align}
Their function does not directly map $n$ to $p(n)$. Instead, the primes $p(n)$ are interspersed with variable-length sequences of $2$s, as determined by the gaps between consecutive primes:
\begin{align*}
\operatorname{Image}(f) = \{ f(n) : n \in \mathbb{N} \} = \{ 2,2,2,3,2,5,2,7,2,2,2,11,2,13,2,\ldots \} .
\end{align*}
Though not technically formulas, Jones et al. and Matiyasevich have constructed various Diophantine equations $D(\vec{x})$ such that the set of prime numbers $\mathbb{P}$ is identical with the set of non-negative values taken on by the polynomial $D(\vec{x})$ as its variables $\vec{x}$ range over the positive integers \cite{jones1976diophantine, matiyasevich1981primediophantineequation}. In the construction of our arithmetic terms for $p(n)$, we obtain several exponential Diophantine equations $E(n,\vec{x}) = 0$ such that, for a given $n \in \mathbb{Z}^+$, the \textit{number of solutions} as the variables $\vec{x}$ range over the natural numbers is equal with $p(n)$ (see Appendix \cref{appendix:F1monomials} for an example). Our equations are the first to be discovered with this remarkable property.

\subsection{Open questions concerning prime formulas}
The challenge to find a formula for the $n$-th prime was finally formalized in the first edition of \textit{An Introduction to the Theory of Numbers}, in which Hardy and Wright posed the following open questions concerning prime numbers:
\begin{question} \label{question:hardyandwright:1} (Hardy and Wright, \cite{hardywright1938introduction})
    Is there a formula for the $n$-th prime?
\end{question}
\begin{question} \label{question:hardyandwright:2} (Hardy and Wright, \cite{hardywright1938introduction})
    Is there a formula for a prime, given the preceding prime?
\end{question}
The construction of an explicit formula for the $n$-th prime is considered fundamentally difficult by various authors. For an information-theoretically approach, see \cite{kolpakov2024impossibility}. While Willans claimed that his formula \cref{eq:willans} addressed Hardy and Wright's questions (one might argue that he was correct, given the phrasing), the fourth edition of \textit{An Introduction to the Theory of Numbers} attempted to clarify the intent and status of the questions, specifying that any formula for $p(n)$ must not utilize any ``\textit{previous knowledge}'' of $p(n)$ \citep[pg.~5]{hardywright1975introduction}. It was further noted that, as of that edition, ``\textit{no satisfactory answer is known}'' \citep[pg.~19]{hardywright1975introduction}.

We propose a constructive solution to \cref{question:hardyandwright:1} with an arithmetic term for $n$-th prime $p(n)$ that does not rely on any previous knowledge of $p(n)$ nor its primality. From such term, we find immediately a recurrence relation to calculate $p(n+1)$ given $p(n)$, thereby providing a constructive solution to \cref{question:hardyandwright:2}. The arithmetic terms we find to solve these problems are immense in size and computationally impractical. However, they are composed of a finite and fixed number of elementary arithmetic operations and universally applicable to all $n$. The status of these questions remains unclear, and it is impossible to know for certain what Hardy and Wright meant by their original questions, however we believe that our results provide satisfactory answers and a definitive resolution. As is often the case, our results also raise new important questions, most notably:
\begin{question}
Does a simpler arithmetic term for the $n$-th prime $p(n)$ exist, or is our formula's great size due to the inherent complexity of the primes themselves?
\end{question}
\begin{question}
Can arithmetic terms for $p(n)$ and $\pi(n)$ be constructed without the hypercube method?
\end{question}
Although our arithmetic terms for $\pi(n)$ and $p(n)$ are large, when one views their great size against the backdrop of the infinitude of the primes, they appear rather small. We hope the study these terms and the search for potential simplifications will lead to new and important discoveries in number theory and mathematics.

\section{Preliminaries}
Under $\mathbb{N}$ we understand the set of natural numbers including $0$.

We introduce the notation $(\vec{x}, \vec{y})$ to denote the \textbf{concatenation} of the tuples $\vec{x}$ and $\vec{y}$. Specifically, if $\vec{x} = (x_1, x_2, \ldots, x_k)$ and $\vec{y} = (y_1, y_2, \ldots, y_j)$, where $k, j \in \mathbb{N}$ represent their respective lengths, then their concatenation is defined as:
\begin{align*}
(\vec{x}, \vec{y}) = (x_1, x_2, \ldots, x_k, y_1, y_2, \ldots, y_j) .
\end{align*}

\subsection{Number theoretic arithmetic terms} \label{SectNumberTheoreticTerms}
The following number theoretic arithmetic terms are used by Mazzanti and Marchenkov in \cite{mazzanti2002plainbases, marchenkov2007superposition}:
\begin{align}
\binom{a}{b} &= \floor{\frac{(2^a+1)^a}{2^{ab}}} \bmod 2^a , \\
\gcd(a,b) &= \floor{\frac{(2^{a^2 b(b+1)} - 2^{a^2 b}) (2^{a^2 b^2} - 1)}{(2^{a^2 b} - 1)(2^{ab^2}-1)2^{a^2 b^2}}} \bmod 2^{ab} , \\
\nu_2(n) &= \floor{\frac{\gcd(n, 2^n)^{n+1} \bmod (2^{n+1}-1)^2}{2^{n+1}-1}} , \\
\HW(n) &= \nu_2\left(\binom{2n}{n}\right) .
\end{align}
Here, $\nu_2(n)$ represents the \textbf{$2$-adic order} of $n$, which is highest exponent of $2$ dividing $n$. $\HW(n)$ denotes the \textbf{Hamming weight} of $n$, which is the number of $1$s in the binary representation of $n$. $\gcd(a,b)$ is the \textbf{greatest common divisor} of $a$ and $b$. A much simpler arithmetic term for $\gcd(a,b)$ was proposed by Prunescu and Shunia in \cite{prunescushunia2024gcd}:
\begin{align}
\gcd(a,b) = \left ( \floor{\frac{2^{ ab(ab + a + b)}}{(2^{a^2 b} - 1)(2^{ab^2}-1)}} \bmod 2^{ab} \right ) - 1 .
\end{align}
The above arithmetic term for the \textbf{binomial coefficient} $\binom{a}{b}$ is well-known and was originally proved by Robinson in \cite{robinson1952arithmetic}. In \cref{subsection:binomialcoeffterm}, we prove two new arithmetic terms for $\binom{a}{b}$ that are of an entirely new construction and greatly simplify our final arithmetic term for $p(n)$.

\subsection{Generalized geometric progressions and the hypercube method} \label{section:hypercube}
Consider $q,r,t \in \mathbb{N}$ such that $q > 1$, $r \geq 0$ and $t \geq 0$. The arithmetic term for the geometric progression
\begin{align*}
\sum_{j=0}^{t-1} q^j = \frac{q^{t}-1}{q-1}
\end{align*}
is well-known. Perhaps lesser known, are the \textbf{generalized geometric progressions of the $r$-th kind}, which are defined as
\begin{align}
G_r (q, t) = \sum_{j=0}^{t-1} j^r q^j .
\end{align}
As described by Matiyasevich in the appendix of \cite{matiyasevich1993hilbert}, for all $r > 0$, $G_r (q, t)$ can be calculated effectively via the recurrence formula:
\begin{align} \label{TermG}
 G_{r} ( q , t ) = \frac{\partial}{\partial q} G_{r-1} ( q , t + 1 ) - \sum_{j=0}^{r-1} \binom{r}{j} G_j ( q , t ) .
\end{align}
Every $G_r(q, t)$ is an arithmetic term in $q$ and $t$.

The \textbf{hypercube method} was discovered by Mazzanti \cite{mazzanti2002plainbases} and is also used in \cite{marchenkov2007superposition, prunescu2024numbertheoreticfunctions}.

Consider $a,b \in \mathbb{N} : 0 \leq a < 2^b$. We define the function
\begin{align}
\delta(a,b) := (2^b - 1)(2^b - a + 1) = 2^{2b} - 2^b a + a - 1.
\end{align}
The Hamming weight of $\delta(a,b)$, denoted by $\HW(\delta(a,b))$, satisfies
\begin{align*}
\HW(\delta(a,b)) = \begin{cases}
2b, & \text{if } a = 0, \\
b, & \text{if } a \neq 0.
\end{cases}
\end{align*}
Let $\vec{n} \in \mathbb{N}^s$ and let $u(\vec{n}),t(\vec{n})$ be arithmetic terms. Now, consider the integer lattice points contained in the $k$-dimensional cube $[0, t(\vec{n})-1]^k$. Define the function
\begin{align*}
  f : [0, t(\vec{n})-1]^k \cap \mathbb{N}^k \to \mathbb{N}
\end{align*}
and assume that
\begin{align*}
    \forall \vec{x} \in [0, t(\vec{n})-1]^k \cap \mathbb{N}^k, \quad f(\vec{x}) < 2^{u(\vec{n})} .
\end{align*}
Define $\beta(\vec{n},\vec{x})$ as the function that maps each point $\vec{x} = (a_1, a_2, \dots, a_k) \in \{0, \dots, t(\vec{n})-1\}^k$ to the integer
\begin{align*}
\beta(\vec{n},\vec{x}) = a_1 + a_2 t(\vec{n}) + \cdots + a_k t(\vec{n})^{k-1}.
\end{align*}
Observe that $\beta(\vec{n},\vec{x})$ provides a bijective mapping between the points in $\{0, \dots, t(\vec{n})-1\}^k$ and the integers from $0$ to $t(\vec{n})^k - 1$. In other words, $\beta(\vec{n},\vec{x})$ enumerates the elements of $\{0, \dots, t(\vec{n})-1\}^k$ in lexicographical order, assigning values from $0$ to $t(\vec{n})^k - 1$.
Let
\begin{align*}
W(\vec{n}) = \sum_{\vec{a} \in \{0, \dots, t(\vec{n})-1\}^k} 2^{2u(\vec{n}) \beta(\vec{n},\vec{a})} \delta(f(\vec{a}), u(\vec{n})).
\end{align*}
We observe that the binary representation of $W(\vec{n})$ corresponds to the concatenation of the binary representations of the numbers $\delta(f(\vec{a}), u(\vec{n}))$ for each $\vec{a} \in \{0, \dots, t(\vec{n})-1\}^k$.

Let $d(\vec{n})$ denote the cardinality of the set $\{ \vec{a} \in \{0, \dots, t(\vec{n}) - 1\}^k : f(\vec{a}) = 0 \}$. It follows that the Hamming weight of $W(\vec{n})$ is given by
\begin{align*}
\HW(W(\vec{n})) = 2 u(\vec{n}) d(\vec{n}) + (t(\vec{n})^k - d(\vec{n})) u(\vec{n}),
\end{align*}
which implies
\begin{align*}
d(\vec{n}) = \frac{\HW(W(\vec{n}))}{u(\vec{n})} - t(\vec{n})^k.
\end{align*}
Therefore, if $W(\vec{n})$ could be expressed as an arithmetic term in $t(\vec{n})$ and $u(\vec{n})$, then the number of zeros of the function $P(\vec{n},\vec{x})$ could also be expressed in such terms. This scenario occurs when $P(\vec{n},\vec{x})$ is an exponential polynomial \textbf{simple in-$\vec{x}$}.

We define a \textbf{simple monomial in-$\vec{x}$} as an expression of the form
\begin{align*}
c v_1^{x_1} \cdots v_k^{x_k} x_1^{r_1} \cdots x_k^{r_k},
\end{align*}
where $r_1, \ldots, r_k \geq 0$, $v_1, \ldots, v_k \geq 1$ are integers, and $c \in \Z$. An \textbf{exponential polynomial simple in-$\vec{x}$} is defined as a sum of such simple monomials.

We apply the identity
\begin{align*}
\sum_{\vec{a} \in \{0, \dots, t(\vec{n})-1\}^k}
a_1^{r_1} v_1^{a_1}
\cdots a_k^{r_k} v_k^{a_k} = G_{r_1}(v_1, t(\vec{n}))
\cdots G_{r_k}(v_k, t(\vec{n}))
= \prod_{i=1}^k G_{r_i}(v_i, t(\vec{n})) ,
\end{align*}
where $G_{r_i}(v_i, t(\vec{n}))$ is the sum function corresponding to each variable $a_i$. The contribution of a simple exponential monomial $m(\vec{n},\vec{x})$ to $W(\vec{n})$ takes the form:
\begin{align} \label{TermA}
\ATerm_k(m(\vec{n},\vec{x}), t(\vec{n}), u(\vec{n}))
&= -(2^{u(\vec{n})} - 1) \cdot c
\cdot
\prod_{i=1}^k G_{r_i}(2^{2u(\vec{n}) t(\vec{n})^{i-1}} v_i, t(\vec{n}))
\end{align}
which is an arithmetic term in $t(\vec{n})$ and $u(\vec{n})$.

If the exponential polynomial contains a free term, meaning $v_1 = \cdots = v_k = 1$ and $r_1 = \cdots = r_k = 0$, the contribution simplifies to
\begin{align} \label{TermC}
\CTerm_k(m(\vec{n},\vec{x}), t(\vec{n}), u(\vec{n}))
= \frac{(2^{u(\vec{n})} - c + 1)(2^{2u(\vec{n}) t(\vec{n})^k} - 1)}{2^{u(\vec{n})} + 1}.
\end{align}
Notice that even for $c = 0$, the contribution of the free term is nonzero. Hence, for exponential polynomials simple in-$\vec{x}$, the quantity $W(\vec{n})$ can always be expressed as an arithmetic term in $t(\vec{n})$ and $u(\vec{n})$.

To summarize, given a parameter tuple $\vec{n} = (n_1, \dots, n_s) \in \mathbb{N}^s$ and a non-negative exponential polynomial function $P(\vec{n},\vec{x})$ simple in-$\vec{x}$, where $\vec{x} = (x_1,\ldots,x_k)$, such that $P(\vec{n},\vec{x})$ is defined on the integer lattice points within the $k$-dimensional cube $[0, t(\vec{n})-1]^k$ and is strictly bounded by $2^{u(\vec{n})}$, the number of zeros within the cube can be expressed as an arithmetic term $d(\vec{n})$ in $t(\vec{n})$ and $u(\vec{n})$.

In practical applications, the coefficients and exponential bases, denoted by $c$ and $v_i$ from the various exponential monomials, will generally depend on some parameter tuple $\vec{n}$ and will also be expressed as arithmetic terms. However, the exponents $r_1, \dots, r_k$ in each monomial are treated as \textbf{constants}. In such cases, appropriate bounds $t(\vec{n})$ and $u(\vec{n})$ are computed so that all relevant zeros lie within the cube $[0, t(\vec{n}) - 1]^k$, and the positive exponential polynomial function remains bounded by $2^{u(\vec{n})}$. Consequently, the number of integer tuples $\vec{x}$ satisfying the equation $P(\vec{n},\vec{x}) = 0$ will be given by an arithmetic term $d(\vec{n})$.

\subsection{Sums of squares of multivariate polynomials}
The construction of an arithmetic term for the $n$-th prime number will require various techniques used in the study of Diophantine equations. To ensure clarity in the subsequent sections, we state a well-known lemma which will be used readily:
\begin{lemma} \label{proof:sumofsquares}
Let $\vec{x} = (x_1,\ldots,x_n)$. Consider polynomials $f_1(\vec{x}), \ldots, f_n(\vec{x}) \in \mathbb{R}[\vec{x}]$, such that for some set $S \subset \mathbb{R}^n$:
\begin{align*}
\forall \vec{x} \in S, \quad
f_1(\vec{x}) \geq 0, \ldots, f_n(\vec{x}) \geq 0 .
\end{align*}
In this case, the sets
\begin{align*}
\{ \vec{x} \in S : f_1(\vec{x}) = 0 \land \ldots \land f_n(\vec{x}) = 0 \}
\end{align*}
and
\begin{align*}
\{ \vec{x} \in S : f_1(\vec{x}) + \cdots + f_n(\vec{x}) = 0 \}
\end{align*}
coincide.
\end{lemma}
We apply this lemma for $f_1(\vec{x}),\ldots,f_n(\vec{x})$ being sums of squares of multivariate polynomials and for $S = \mathbb{N}^n$.

\subsection{Exponential Diophantine single-fold definitions} \label{SectDef}

Every of the following definitions will be the square of an exponential polynomial expression, or a sum of squares of such expressions. We use the notation:
\begin{align*}
    E(\vec x, [k]) = 0,
\end{align*}
to express the fact that $k$ many variables are quantified existentially. We can express conjunctions of such definitions in the form:
\begin{align*}
E(\vec x, [k]) + E(\vec y, [m]) = 0,
\end{align*}
and it is always tacitly understood, that none of the $k$ many quantified variables from the first definition appear in the set of $m$ quantified variables of the second definition. The number of quantified variables will be $k + m$. There are however situations in which one combine definitions and it is important that some variable appears in two or more such expressions. In this case, these variables has to be written explicitly.

A relation $R(\vec x, y)$ has a {\bf single-fold} (exponential) Diophantine definition
\begin{align*}
    E(\vec{x}, [k], y) = 0 ,
\end{align*}
where $[k]$ means a tuple of variables $\vec y = (y_1, \dots, y_k)$, if and only if:
\begin{enumerate}
    \item[(i)] For all $\vec x \in \mathbb N^n$ and for all $y \in \mathbb N$, $R(\vec x, y)$ is true if and only if
    $$\exists! \vec y : E(\vec{x}, \vec y, y) = 0.$$
    \item[(ii)] If $R(\vec x, y)$ is true, then the corresponding tuple $\vec y$ satisfying this equation is uniquely determined by the tuple $(\vec x, y)$.
\end{enumerate}
In the special situation that the relation $R(\vec x, y)$ is a function $f(\vec x) = y$, the value of $y$ is  uniquely determined by the value of $\vec x$ as well. Consequently, in this situation the whole tuple $(\vec y, y)$ is uniquely determined by the value of the tuple $\vec x$. We will deal with several single-fold (exponential) Diophantine definitions and most of them define functions.

Now, we will define various single-fold relations with corresponding lemmas. We omit most of the lemma proofs, since they are quite trivial and follow from their definitions and preceding statements.

We always need the operations quotient and remainder. The single-fold relation $z = \floor{x / y}$ is denoted as
\begin{align*}
E_/ (x, y, [2], z) = 0 .
\end{align*}
\begin{lemma} \label{proof:divsingle-fold}
$
\forall (x_1, x_2, x_3) \in \mathbb{N}^3, \quad
x_3 = \floor{x_1 / x_2}
\iff E_/ (x_1, x_2, [2], x_3) = 0
$
\begin{align*}
\iff \exists \vec{y} \in \mathbb{N}^2 : (x_1 - x_2 x_3 - y_1)^2 + (y_1 + y_2 + 1 - x_2)^2 = 0 .
\end{align*}
\end{lemma}

The single-fold relation $z = x \bmod y$ is denoted as
\begin{align*}
E_{\bmod} (x, y, [2], z) = 0 .
\end{align*}
\begin{lemma} \label{proof:modsingle-fold}
$
\forall (x_1, x_2, x_3) \in \mathbb{N}^3, \quad
x_3 = x_1 \bmod x_2
\iff E_{\bmod}(x_1, x_2, [2], x_3) = 0
$
\begin{align*}
\iff \exists \vec{y} \in \mathbb{N}^2 : (x_1 - x_2 y_1 - x_3)^2 + (x_3 + y_2 + 1 - x_2)^2 = 0 .
\end{align*}
\end{lemma}

A single-fold relation of the divisibility condition $y \mid x$ will be denoted with
\begin{align*}
E_{|}(x, [1], y) = 0 .
\end{align*}
\begin{lemma}
$
\forall (x_1, x_2) \in \mathbb{N}^2, \quad
x_2 \mid x_1
\iff E_{|}(x_1, [1], x_2) = 0
\iff \exists y_1 \in \mathbb{N} : (x_1 - x_2 y_1)^2 = 0 .
$
\end{lemma}

A single-fold relation of the indivisibility condition $y \nmid x$ will be denoted with
\begin{align*}
E_{\nmid}(x, [3], y) = 0 .
\end{align*}
\begin{lemma}
$
\forall (x_1, x_2) \in \mathbb{N}^2, \quad
x_2 \nmid x_1
\iff E_{\nmid}(x_1, [3], x_2) = 0
$
\begin{align*}
\iff \exists \vec{y} \in \mathbb{N}^3 : (x_1 - x_2 y_1 - y_2 - 1)^2 + (y_2 + y_3 + 2 - x_2)^2 = 0 .
\end{align*}
\end{lemma}

The expression $y = \nu_2(x)$ means that $y$ is the exponent of $2$ in the prime number decomposition of $x$. We will denote this single-fold relation with
\begin{align*}
E_{\nu}(x, [4], y) = 0 .
\end{align*}
\begin{lemma} \label{proof:nusingle-fold}
$
\forall (x_1, x_2) \in \mathbb{N}^2, \quad
x_2 = \nu_2(x_1)
\iff E_{\nu}(x_1, [4], x_2) = 0
$
\begin{align*}
\iff \exists \vec{y} \in \mathbb{N}^4 : E_\nmid (x_1, [3], 2^{x_2+1}) + E_| (x_1, [1], 2^{x_2}) = 0 .
\end{align*}
\end{lemma}

A requirement of the hypercube method (\cref{section:hypercube}) is that all exponentiations within relations must be powers of a fixed integer base (most commonly $2$) with an exponent that is \textbf{simple in all variables}, meaning that the exponent is a linear function in the displayed unknowns \cite{mazzanti2002plainbases, marchenkov2007superposition}. To illustrate, $2^{2 x + 3 y}$ is a valid exponentiation for a hypercube relation, while $x^{y}$ is invalid. Therefore, we will require a general relation which translates exponentiations of the form $z = x^y$ to an equivalent form that is suitable for the hypercube method. The single-fold relation $z=x^y$ will be written as
\begin{align*}
E_{\exp}(x, y, [4], z) = 0 .
\end{align*}

\begin{lemma}
$ \forall (x_1, x_2, x_3) \in \mathbb{N}^3, \quad
x_3 = x_1^{x_2}
\iff E_{\exp}(x_1, x_2, [4], x_3) = 0 $
\begin{align*}
\iff \exists
\vec{y} \in \mathbb{N}^4 : (y_1 - x_1 x_2 - x_1 - 1)^2 + (y_2 - y_1 x_2)^2 + E_{\bmod} (2^{y_2}, 2^{y_1} - x_1, [2], x_3) = 0 .
\end{align*}
\end{lemma}
\begin{proof}
From Mazzanti \cite{mazzanti2002plainbases}, the general exponentiation can be computed using only powers of $2$ by the formula:
\begin{align*}
x_1^{x_2} = 2^{(x_1 x_2 + x_1 + 1) x_2} \bmod (2^{x_1 x_2 + x_1 + 1} - x_1) .
\end{align*}
This means that
\begin{align*}
x_3 = x_1^{x_2} \iff E_{\bmod} ( 2^{(x_1 x_2 + x_1 + 1) x_2} , 2^{x_1 x_2 + x_1 + 1} - x_1, [2], x_3) = 0 .
\end{align*}
To ensure that our definitions are simple in all variables, meaning that the powers of $2$ are linear functions in the displayed unknowns, we introduce a new variables $y_1 = x_1 x_2 + x_1 + 1$, $y_2 = y_1 x_2$ as the sum of squares $(y_1 - x_1 x_2 - x_1 - 1)^2 + (y_2 - y_1 x_2)^2 = 0 $.
After re-writing the exponents for $2^{(x_1 x_2 + x_1 + 1) x_2}$ and $2^{x_1 x_2 + x_1 + 1} = 2^{y_1}$ in terms of $y_1,y_2$, we obtain
\begin{align*}
E_{\exp}(x_1, x_2, [4], x_3)
= (y_1 - x_1 x_2 - x_1 - 1)^2 + (y_2 - y_1 x_2)^2 + E_{\bmod} (2^{y_2}, 2^{y_1} - x_1, [2], x_3) = 0 ,
\end{align*}
which defines $x_3 = x_1^{x_2}$ single-fold.
\end{proof}
\section{Single-fold definitions for binomial coefficients} \label{section:binomialcoefficients}
We will need a single-fold exponential Diophantine definition of the relation $z = \binom{x}{y}$ for two important constructions: The single-fold Diophantine definitions for the factorial function and respectively for the Hamming weight of a natural number. For our initial construction, we will use the fact that \cite{robinson1952arithmetic}:
\begin{align*}
z = \binom{x}{y} = \floor{\frac{(2^x+1)^x}{2^{xy}}} \bmod 2^x .
\end{align*}
We denote this version of the single-fold relation $z = \binom{x}{y}$ as
\begin{align*}
    E_{\binom{\#}{\#}}(x, y, [12], z) = 0 .
\end{align*}

\begin{lemma}
$
\forall (x_1, x_2, x_3) \in \mathbb{N}^3, \quad
x_3 = \binom{x_1}{x_2}
\iff E_{\binom{\#}{\#}}(x_1, x_2, [12], x_3) = 0 $
\begin{align*}
\iff \exists \vec{y} \in \mathbb{N}^{12} : \quad
& (1+2^{x_1}-y_1)^2
+ E_{\exp}(y_1, x_1, [4], y_2)
+ (y_3 - x_1 x_2 )^2
+ E_{/}(y_2, 2^{y_3}, [2], y_4) \\
& + E_{\bmod}(y_4, 2^{x_1}, [2], x_3) = 0 .
\end{align*}
\end{lemma}
This version contains $12$ quantified variables. Reducing the number of quantified variables can vastly decrease the magnitude of arithmetic terms constructed using the hypercube method. The single-fold relation $z = \binom{x}{y}$ is important and will be used in the definition of additional single-folds, so it would be ideal if we could see some reduction here.

Indeed, we find a new arithmetic term for $z = \binom{x}{y}$ that allows us to reduce the number of quantified variables from $12$ to $7$. This new version will be written as
\begin{align*}
E_{\binom{\#}{\#}}(x,y,[7],z) = 0
\end{align*}
and will be derived in \cref{subsection:binomialcoeffterm} below.

\subsection{A new arithmetic term for binomial coefficients} \label{subsection:binomialcoeffterm}
Consider the \textbf{Padovan sequence}, whose terms are the integers given by the recurrence relation
\begin{align*}
P(k) = P(k-2) + P(k-3) ,
\end{align*}
with initial starting conditions $P(0)=1$, $P(1)=P(2)=0$, $P(3)=1$. For the sequence terms, see \seqnum{A000931} in the OEIS.

Let us now generalize this type of recurrence relation as
\begin{align*}
s_d(n) = s_d(n-d+1) + s_d(n-d) ,
\end{align*}
with initial starting conditions $s_d(0)=s_d(1)=\cdots=s_d(d-2)=0$, $s_d(d-1)=1$. We call $s_d(n)$ a \textbf{generalized Padovan sequence}. From this definition, $s_3(k) = P(k+1)$ for all $k \in \mathbb{N}$.

We represent the finite sequence of the first $d^2$ elements of the sequence of degree $d$, in the following way: The sequence is cut in $d$ many segments and they are arranged in $d \times d$ matrix.

For example, for $d=8$, one gets:
\begin{align*}
\begin{array}{llllllllll}
s_8(n) =
& 0, & 0, & 0, & 0, & 0, & 0, & 0, & 1, \\
& 0, & 0, & 0, & 0, & 0, & 0, & 1, & 1, \\
& 0, & 0, & 0, & 0, & 0, & 1, & 2, & 1, \\
& 0, & 0, & 0, & 0, & 1, & 3, & 3, & 1, \\
& 0, & 0, & 0, & 1, & 4, & 6, & 4, & 1, \\
& 0, & 0, & 1, & 5, & 10, & 10, & 5, & 1, \\
& 0, & 1, & 6, & 15, & 20, & 15, & 6, & 1, \\
& 1, & 7, & 21, & 35, & 35, & 21, & 7, & \textcolor{red}{2}, \\
& \vdots & \vdots & \vdots & \vdots & \vdots & \vdots & \vdots & \vdots
\end{array}
\end{align*}
Observe that Pascal's triangle is embedded in this matrix, as
the recurrence rule of the sequence simulates the recurrence rule
of the binomial coefficients:
\begin{align*}
\binom{a}{b} = \binom{a-1}{b-1} + \binom{a-1}{b} .
\end{align*}
So, the row number $7$ of the matrix contains the row number $6$ of Pascal's Triangle, while the $64$-th element $s_8(63)$ is equal to $2$, as from now on the row-sides collide and the simulation ends. To illustrate:
\begin{align*}
\begin{array}{llllllllll}
s_8(n) =
& 0, & 0, & 0, & 0, & 0, & 0, & 0, & \binom{0}{0}, \\
& 0, & 0, & 0, & 0, & 0, & 0, & \binom{1}{0}, & \binom{1}{1}, \\
& 0, & 0, & 0, & 0, & 0, & \binom{2}{0}, & \binom{2}{1}, & \binom{2}{2}, \\
& 0, & 0, & 0, & 0, & \binom{3}{0}, & \binom{3}{1}, & \binom{3}{2}, & \binom{3}{3}, \\
& 0, & 0, & 0, & \binom{4}{0}, & \binom{4}{1}, & \binom{4}{2}, & \binom{4}{3}, & \binom{4}{4}, \\
& 0, & 0, & \binom{5}{0}, & \binom{5}{1}, & \binom{5}{2}, & \binom{5}{3}, & \binom{5}{4}, & \binom{5}{5}, \\
& 0, & \binom{6}{0}, & \binom{6}{1}, & \binom{6}{2}, & \binom{6}{3}, & \binom{6}{4}, & \binom{6}{5}, & \binom{6}{6}, \\
& \binom{7}{0}, & \binom{7}{1}, & \binom{7}{2}, & \binom{7}{3}, & \binom{7}{4}, & \binom{7}{5}, & \binom{7}{6}, & \textcolor{red}{2}, \\
& \vdots & \vdots & \vdots & \vdots & \vdots & \vdots & \vdots & \vdots
\end{array}
\end{align*}

In order to produce an arithmetic term able to represent the row number $d-1$
of this table, we consider following polynomials:
\begin{align*}
S_{d,k}(x) &= s_d(0) x^k + s_d(1) x^{k-1} + \cdots + s_d(0) , \\
B_d(x) &= x^d - x - 1 .
\end{align*}

\begin{theorem} \label{proof:binomialcoeffmodmod}
\begin{align*}
\forall a,b \in \mathbb{N}, \quad
\binom{a}{b} = \left( 2^{2(a+2)((a+1)^2+b+1)} \bmod (2^{2(a+2)^2}-2^{2(a+2)}-1) \right) \bmod 2^{2(a+2)} .
\end{align*}
\end{theorem}
\begin{proof}
Under the conditions $s_d(0) = s_d(1) = \cdots = s_d(d-2) = 0$ and $s_d(d-1) = 1$,
and taking $k \geq d$, one has that:
\begin{align*}
B_d(x) S_{d,k}(x) &= x^{k+1} - (s_d(k-d+2) + s_d(k-d+1)) x^{d-1} \notag \\
& - (s_d(k-d+3) + s_d(k-d+2)) x^{d-2} \notag \\
& - \cdots - (s_d(k) + s_d(k-1)) x - s_d(k) .
\end{align*}

Applying the recurrence rule for the sums in parentheses, we find:
\begin{align*}
& s(k-d+2) + s(k-d+1) = s(k+1), \\
& s(k-d+3) + s(k-d+2) = s(k+2), \\
& \vdots \\
& s(k) + s(k-1) = s(k + d - 1),
\end{align*}
so
\begin{align*}
B_d(x) S_{d,k}(x) = x^{k+1} - s_d(k+1) x^{d-1} - \cdots - s_d(k+d-1) x - s_d(k) .
\end{align*}

For every value $x \in \mathbb N$, one gets the congruence:
\begin{align*}
x^{k+1} \equiv s_d(k+1) x^{d-1} + \cdots + s_d(k+d-1) x + s_d(k) \pmod{x^d - x - 1} .
\end{align*}

We are looking for a value of $x$ such that
\begin{align*}
0 \leq s_d(k+1) x^{d-1} + \cdots + s_d(k+d-1) x + s_d(k) < x^d - x - 1 ,
\end{align*}
because we want that
\begin{align*}
x^{k+1} \bmod (x^d-x-1) = s_d(k+1) x^{d-1} + \cdots + s_d(k+d-1)x + s_d(k) .
\end{align*}
We observe that for $k$ in the row $d-1$, the elements $s_d(k+1), \dots, s_d(k+d-1)$
belong to either row number $d-1$ or to row number $d$, so they are all $< 2^{d-1}$. We take $x = 4^d$.

Indeed, in order to prove
\begin{align*}
s_d(k+1) (4^d)^{d-1} + \cdots + s_d(k+d-1) 4^d + s_d(k) + 4^d + 1 < 4^{d^2}
\end{align*}
it is enough to show that
\begin{align*}
2^{d-1} (4^d)^{d-1} + \cdots + 2^{d-1} 4^d + 2^{d-1} + 4^d + 1 < 4^{d^2} .
 \end{align*}
The biggest term (leftmost) is
\begin{align*}
2^{d-1} (4^d)^{d-1} = 2^{2d^2 - 2d + d - 1} = 2^{2d^2 - d - 1}
\end{align*}
and there are exactly $d+1$ terms, so it suffices to show that
\begin{align*}
& (d + 1) 2^{2d^2 - d - 1} < 2^{2d^2} , \\
\Leftrightarrow
& \log_2(d+1) + 2d^2 - d - 1 < 2d^2 , \\
\Leftrightarrow
& \log_2(d+1) < 2d^2 ,
\end{align*}
which is true for all integers $d \geq 0$.

Thus, we find out that
\begin{align*}
4^{d(k+1)} \bmod (4^{d^2} - 4^d - 1) = s_d(k+1) (4^d)^{d-1} + \cdots + s_d(k+d-1) 4^d + s_d(k) ,
\end{align*}
hence
\begin{align*}
s_d(k) = \left( 4^{d(k+1)} \bmod (4^{d^2}-4^d-1) \right) \bmod 4^d.
\end{align*}
for $d \leq k < d^2-d-1$.

In order to compute $\binom{a}{b}$ with $0 \leq b \leq a$, we take
$d=a+2$ and $k=(a+1)^2+b$ and we find that
\begin{align*}
\binom{a}{b} &= \left( 4^{(a+2)((a+1)^2+b+1)} \bmod (4^{(a+2)^2}-4^{a+2}-1) \right) \bmod 4^{a+2} \\
&= \left( 2^{2(a+2)((a+1)^2+b+1)} \bmod (2^{2(a+2)^2}-2^{2(a+2)}-1) \right) \bmod 2^{2(a+2)} .
\end{align*}
\end{proof}

\begin{theorem} \label{proof:binomialcoeffdivmod}
\begin{align*}
\forall a,b \in \mathbb{N}, \quad
\binom{a}{b} = \floor{\frac{2^{2(a+2)((a+1)^2+b+1)}}{2^{2(a+2)^2} - 2^{2(a+2)} - 1}} \bmod 2^{2(a+2)} .
\end{align*}
\end{theorem}
\begin{proof}
We start with the formula proved in \cref{proof:binomialcoeffmodmod}:
\begin{align*}
\binom{a}{b} = \left( 2^{2(a+2)((a+1)^2+b+1)} \bmod (2^{2(a+2)^2}-2^{2(a+2)}-1) \right) \bmod 2^{2(a+2)} .
\end{align*}

Now, let $x = 2^{2(a+2)((a+1)^2+b+1)}$ and let $y = 2^{2(a+2)^2}-2^{2(a+2)}-1$.

Substituting in the well-known identity $x \bmod y = x - y \floor{x/y}$, we obtain
\begin{align*}
\binom{a}{b} = ( x \bmod y ) \bmod 2^{2(a+2)} = \left( x - y \floor{\frac{x}{y}} \right) \bmod 2^{2(a+2)} .
\end{align*}
As the application of $\bmod$ defined on the ring of integers $\Z$ and with values in the ring of remainder classes $\mathbb \Z / 2^{2(a+2)} \Z$, is a homomorphism of rings, we see that
\begin{align*}
\binom{a}{b}
&= \left (x \bmod 2^{2(a+2)}\right ) - \left (y \bmod 2^{2(a+2)} \right ) \left (\floor{\frac{x}{y}} \bmod 2^{2(a+2)}\right ) \\
&= (0) - (-1) \floor{\frac{x}{y}} \bmod 2^{2(a+2)} \\
&= \floor{\frac{x}{y}} \bmod 2^{2(a+2)} .
\end{align*}
Substituting back $x$ and $y$, we have
\begin{align*}
\binom{a}{b} &= \floor{\frac{2^{2(a+2)((a+1)^2+b+1)}}{2^{2(a+2)^2} - 2^{2(a+2)} - 1}} \bmod 2^{2(a+2)} .
\end{align*}
\end{proof}

\subsection{An improved single-fold definition for binomial coefficients}

Since all exponentiations in the arithmetic terms from \cref{proof:binomialcoeffmodmod} and \cref{proof:binomialcoeffdivmod} are performed using powers of $2$, we no longer require the single-fold relation $E_{\exp}(x,y,[4],z) = 0$ in our construction of the single-fold for $z = \binom{x}{y}$. This results in an immediate reduction of $4$ quantified variables. We get a further reduction of $1$ quantified variable, for a total reduction of $5$ quantified variables:

\begin{lemma} \label{proof:binomialsingle-fold2}
$
\forall (x_1, x_2, x_3) \in \mathbb{N}^3, \quad
x_3 = \binom{x_1}{x_2}
\iff E_{\binom{\#}{\#}}(x_1, x_2, [7], x_3) = 0
$
\begin{align*}
\iff \exists \vec{y} \in \mathbb{N}^7 : \quad
& (y_1 - (2 x_1^3 + 8 x_1^2 + 2 x_1 x_2 + 12 x_1 + 4 x_2 + 8 ))^2 \\
& + (y_2 - (2 x_1^2 + 8 x_1 + 8))^2
+ E_{/} (2^{y_1}, 2^{y_2} - 2^{2 x_1 + 4} - 1, [2], y_3) \\
& + E_{\bmod} (y_3, 2^{2 x_1 + 4}, [2] , x_3) = 0 .
\end{align*}
\end{lemma}
\begin{proof}
The proof follows immediately from \cref{proof:binomialcoeffdivmod}.
\end{proof}

\begin{remark}
In \cref{proof:binomialsingle-fold2}, we have elected to use the \textbf{div-mod} representation of $\binom{a}{b}$ (\cref{proof:binomialcoeffdivmod}) as opposed to the \textbf{mod-mod} representation (\cref{proof:binomialcoeffmodmod}) because the div-mod version results in fewer monomials after expanding all squares in the equation.
\end{remark}

\section{Single-fold definitions for factorials and Hamming weights} \label{section:factorialsandhw}
We will now develop our single-fold relations for the factorial and Hamming weight functions. We start with some lemmas:

\begin{lemma} \label{proof:factorialr}
$\forall n \in \mathbb{N}, \quad 8^{n^2} \geq (n+1)^{n+2}$ .
\end{lemma}
\begin{proof}
 For $n = 0$ we have $1 \geq 1$. For $n = 1$ we have $8 \geq 8$. For $n = 2$, we have $8^4 \geq 3^4$. Further, the functions are more and more apart, because for $n = 3$, one has already $n^2 > (n+2) \log_8 (n+1)$, and $n^2$ increases faster than the right-hand side.
\end{proof}

\begin{lemma} \label{proof:factorialterm}
\begin{align*}
\forall n \in \mathbb{N}, \quad
n! = \floor
{
2^{3n^3}
\Big/
\left(
\floor{
\frac{2^{2(2^{3n^2}+2)((2^{3n^2}+1)^2+n+1)}}
{2^{2(2^{3n^2}+2)^2} - 2^{2(2^{3n^2}+2)} - 1}}
\bmod 2^{2(2^{3n^2}+2)}
\right)
}
\end{align*}
\end{lemma}
\begin{proof}
From Robinson \cite{robinson1952arithmetic}, we have
\begin{align*}
n! = \lim_{a\to\infty} a^n \Big/ \binom{a}{n}
\end{align*}
and if $a > (2n)^{n+1}$, then
\begin{align*}
n! = \floor{ a^n \Big/ \binom{a}{n} } .
\end{align*}
Matiyasevich gave a slightly different proof in \citep[pg. 46]{matiyasevich1993hilbert}, which is valid for all $a \geq (n+1)^{n+2}$. We opt to use Matiyasevich's bound, since $(2n)^{n+1} \gg (n+1)^{n+2}$ as $n \to \infty$.

Put $a = 8^{n^2}$. By \cref{proof:factorialr}, we have that $\forall n \in \mathbb{N}, \; a = 8^{n^2} \geq (n+1)^{n+2}$.
It follows that
\begin{align*}
n! &= \floor{ (8^{n^2})^n \Big/ \binom{8^{n^2}}{n} } .
\end{align*}
Substituting in the formula for $\binom{8^{n^2}}{n}$ from \cref{proof:binomialcoeffdivmod}, one obtains:
\begin{align*}
n! &= \floor
{
8^{n^3} \Big/
\left(
\floor{\frac{4^{(8^{n^2}+2)((8^{n^2}+1)^2+n+1)}}{4^{(8^{n^2}+2)^2} - 4^{8^{n^2}+2} - 1}}
\bmod 4^{8^{n^2}+2}
\right)
} \\
&= \floor
{
2^{3n^3}
\Big/
\left(
\floor{
\frac{2^{2(2^{3n^2}+2)((2^{3n^2}+1)^2+n+1)}}
{2^{2(2^{3n^2}+2)^2} - 2^{2(2^{3n^2}+2)} - 1}}
\bmod 2^{2(2^{3n^2}+2)}
\right)
}.
\end{align*}
\end{proof}

We define the single-fold exponential Diophantine relation $f = n!$ as:
\begin{align*}
    E_{!}(n, [13], f) = 0 .
\end{align*}

\begin{lemma} \label{proof:factorialsingle-fold}
$
\forall (x_1, x_2) \in \mathbb{N}^2, \quad
x_2 = x_1!
\iff
E_{!}(x_1, [13], x_2) = 0
$
\begin{align*}
\iff \exists \vec{y} \in \mathbb{N}^{13} : \;
& (y_1 - x_1^2)^2
+ (y_2 - 2^{3 y_1})^2
+ (y_3 - x_1 y_1)^2
+ E_{\binom{\#}{\#}}(y_2, x_1, [7], y_4) \\
&+ E_{/}(2^{3 y_3}, y_4, [2], x_2)
= 0 .
\end{align*}
\end{lemma}
\begin{proof}
The proof follows immediately from \cref{proof:factorialterm}.
\end{proof}

A further application of the binomial coefficient is the single-fold definition of the Hamming Weight of a number $y = \HW(x)$. We write this expression as
\begin{align*}
    E_{\HW}(x, [12], y) = 0 .
\end{align*}

\begin{lemma}
$ \forall (x_1, x_2) \in \mathbb{N}^2, \quad
x_2 = \HW(x_1)
\iff
E_{\HW}(x_1, [12], x_2) = 0$
\begin{align*}
\iff \exists \vec{y} \in \mathbb{N}^{12} : E_{\binom{\#}{\#}}(2 x_1, x_1, [7], y_1) + E_\nu(y_1, [4], x_2) = 0 .
\end{align*}
\end{lemma}
\begin{proof}
The lemma follows immediately from the fact that: $\HW(x_1) = \nu_2 \left( \binom{2x_1}{x_1} \right)$.
\end{proof}

\section{The prime omega function and modular square roots of unity} \label{section:omega}
The \textbf{prime omega function}, $\omega(n)$, is defined for positive integers $n$ and returns the number of distinct primes dividing $n$.

The \textbf{modular square roots of unity counting function}, $\NTerm(n)$, is defined for all natural numbers $n$ and counts the number of square roots of unity modulo $n$. For the special case $n=0$, we define $\NTerm(0) = 0$. For all $n > 0$, we define this function as:
\begin{align}
\NTerm(n) = | \{ a \in \{0,\ldots,n-1\} : a^2 \equiv 1 \pmod{n} \} | .
\end{align}
 We note that $\NTerm(1) = 1$, since $0 \equiv 1 \pmod{1}$.

As our first step towards finding an arithmetic term for $\omega(n)$, we will apply the hypercube method to find an arithmetic term for $\NTerm(n)$. This will require the very important arithmetic term:
\begin{align} \label{TermM}
\begin{array}{llll}
M(n)
&= \CTerm_2(1, t(n), u(n))
&+ \ATerm_2(x_1^4, t(n), u(n))
&+ \ATerm_2(-2 x_1^2, t(n), u(n)) \\
&+ \ATerm_2(-2n x_1^2 x_2, t(n), u(n))
&+ \ATerm_2(n^2 x_2^2, t(n), u(n))
&+ \ATerm_2(2n x_2, t(n), u(n)) ,
\end{array}
\end{align}
where $t(n) = n + 1$ and $u(n) = n + 5$. For the full arithmetic term representation of $M(n)$, see \cref{TermMIntegerRepresentation}.

\begin{lemma} \label{proof:ntermarithmetic}
For all $n \in \mathbb{N}$, the function $\NTerm(n)$ is given by the arithmetic term:
\begin{align} \label{TermNArithmetic}
\NTerm(n) = \frac{\HW(M(4n))}{u(4n)} - t(4n)^2 ,
\end{align}
where $t(n)=n+1$, $u(n)=n+5$, and $M(n)$ is the arithmetic term defined in \cref{TermM}.
\end{lemma}
\begin{proof}
Let $n \in \Z^+$. Consider the congruence:
\begin{align*}
a^2 \equiv 1 \pmod{n} .
\end{align*}
In order to find an arithmetic term that expresses the number of solutions to this congruence, we consider the Diophantine equation:
\begin{align*}
x_1^2 - n x_2 = 1 .
\end{align*}
We observe that if $(x_1, x_2)$ is a solution with $x_1 < n$, then
\begin{align*}
    x_2 = \frac{x_1^2 - 1}{n} < \frac{n^2 - 1}{n} < \frac{n^2}{n} = n .
\end{align*}
Clearly, the number of integer pairs $(x_1,x_2)$ satisfying the above equation with $0 \leq x_1, x_2 < n$ equals the number of elements $a \in \{ 0, \ldots, n-1 \}$ such that $a^2 \equiv 1 \pmod{n}$. More formally,
\begin{align*}
| \{ (x_1,x_2) \in \{0,\ldots,n-1\}^2 : x_1^2 -  n x_2 = 1  \} |
\; = \;
| \{ a \in \{0,\ldots,n-1\} : a^2 \equiv 1 \pmod{n} \} | .
\end{align*}
Applying the hypercube method as described in \cref{section:hypercube}, we can count the number of solutions as an arithmetic term.

One must count the number of pairs $(x_1, x_2) \in \{0, 1, \dots, n-1\} \times \{0, 1, \dots, n-1\}$ which satisfy the equation:
\begin{align*}
    (x_1^2 - n x_2 - 1)^2 = 0 ,
\end{align*}
which develops to:
\begin{align*}
    x_1^2 + n^2 x_2^2 + 1 - 2n x_1^2 x_2 - 2 x_1^2 + 2 n x_2 = 0 .
\end{align*}
One can take $t(n) = n + 1$ and $u(n) = n + 5$. We verify by numerical calculations that
\begin{align*}
\forall n \in \Z^+ : n \leq 32, \quad
| \{ (x_1,x_2) \in \{0,\ldots,n-1\}^2 : x_1^2 - n x_2 - 1 < 2^{n+5} \} | = 0 .
\end{align*}
Observe that the distance between the exponential function and the polynomial one is permanently increasing. This tells us
\begin{align*}
\forall n \in \Z^+ , \quad
| \{ (x_1,x_2) \in \{0,\ldots,n-1\}^2 : x_1^2 - n x_2 - 1 < 2^{n+5} \} | = 0 .
\end{align*}
Hence
\begin{align*}
\NTerm(n) = | \{ (x_1,x_2) \in \{0,\ldots,n-1\}^2 : x_1^2 - n x_2 = 1 \} |
= \frac{\HW(M(n))}{n+5} - (n+1)^2
= \frac{\HW(M(n))}{u(n)} - t(n)^2 .
\end{align*}
Finally, by numerical calculations we verify the special case $n = 0$ to find that
\begin{align*}
\frac{\HW(M(0))}{u(0)} - t(0)^2 = \NTerm(0) = 0 .
\end{align*}
Thus, the given arithmetic term is valid for all $n \in \mathbb{N}$.
\end{proof}

\begin{remark}
We could have used $t(n) = 4n$ and $u(n)=4n+4$ in \cref{proof:ntermarithmetic} to obtain an arithmetic term for $\NTerm(n)$ that is valid for all $n \in \Z^+$. However, we opted instead to use $t(n)=4n+1$ and $u(n)=4n+5$ because these values lead to fewer monomials in \cref{section:monomials}.
\end{remark}

\begin{lemma} \label{proof:omegasolutions}
$
\forall n \in \Z^+, \quad
\omega(n) = \nu_2(\NTerm(4n)) - 1 .
$
\end{lemma}
\begin{proof}
Let the prime number decomposition of $n$ be $n = 2^{a} p_1^{b_1} \dots p_d^{b_d}$ where $p_i$ are odd primes, $b_i \geq 1$, $d \geq 0$ and $a \geq 0$.
By the Chinese Remainder Theorem, there is the following isomorphism of rings:
\begin{align*}
\Z/(4n)\Z \cong \Z/(2^{a+2}) \Z \times \Z/(p_1^{b_1})\Z \times \cdots \times \Z/(p_d^{b_d})\Z .
\end{align*}

We use the following known facts:
\begin{enumerate}
\item[(i)] In $\Z/2\Z$, the congruence $x^2 \equiv 1 \pmod{2}$ has exactly $1$ solution:
\begin{align*}
    x = 1 .
\end{align*}
\item[(ii)] In $\Z/4\Z$, the congruence $x^2 \equiv 1 \pmod{4}$ has exactly $2$ solutions:
\begin{align*}
    x = 1, \quad x = 3 .
\end{align*}
\item[(iii)] In $\Z/2^k\Z$ for $k \geq 3$, the congruence $x^2 \equiv 1 \pmod{2^k}$ has exactly $4$ solutions:
\begin{align*}
    x = 1, \quad x = 2^{k-1} - 1, \quad x = 2^{k-1} + 1, \quad x = 2^k - 1 .
\end{align*}
\item[(iv)] In $\Z/p^k\Z$ for $p$ an odd prime and $k \geq 1$, the congruence $x^2 \equiv 1 \pmod{p^k}$ has exactly $2$ solutions:
\begin{align*}
    x = 1, \quad x = p^k - 1 .
\end{align*}
\end{enumerate}

Now, every solution to the equation $x^2 \equiv 1 \pmod{4n}$ corresponds via the Chinese Remainder Theorem isomorphism to a tuple of solutions in the given rings. We proceed with a case discussion.

\textbf{Case 1:} The number $n$ is odd.

In this case, $n = p_1^{b_1} \cdots p_d^{b_d}$ and $\omega(n) = d$. The number of solutions in $\Z/(4n)\Z$ will be:
\begin{align*}
(2 \text{ solutions in } \Z/4\Z) \times (2 \text{ solutions in } \Z/p_1^{b_1}\Z) \times \cdots \times (2 \text{ solutions in } \Z/p_d^{b_d}\Z) = 2^{d+1}.
\end{align*}
Thus, $\omega(n) = d = \nu_2(\NTerm(4n)) - 1$.

\textbf{Case 2:} The number $n$ is even.

In this case, $n = 2^{a} p_1^{b_1} \cdots p_d^{b_d}$, $a \geq 1$ and $\omega(n) = d + 1$. The exponent of $2$ in $4n$ will be $ a + 2 \geq 3$, so the number of solutions in $\Z/(4n)\Z$ will be:
\begin{align*}
(4 \text{ solutions in } \Z/2^{a + 2}\Z) \times (2 \text{ solutions in } \Z/p_1^{b_1}\Z) \times \cdots \times (2 \text{ solutions in } \Z/p_d^{b_d}\Z) = 2^{d+2}.
\end{align*}
Thus, $\omega(n) = d + 1 = \nu_2(\NTerm(4n)) - 1$ again.
\end{proof}

\begin{theorem} \label{proof:omega}
For all $n \in \Z^+$, the number of distinct prime divisors of $n$, $\omega(n)$, is given by the arithmetic term:
\begin{align} \label{OmegaTerm}
\omega(n) = \nu_2\left( \frac{\HW(M(4n))}{u(4n)} - t(4n)^2 \right) - 1 ,
\end{align}
where $t(n)=n+1$, $u(n)=n+5$, and $M(n)$ is the arithmetic term defined in \cref{TermM}.
\end{theorem}
\begin{proof}
Let $n \in \Z^+$. By \cref{proof:omegasolutions} we have $\omega(n) = \nu_2(\NTerm(4n))-1$. Applying \cref{proof:ntermarithmetic}, we obtain the arithmetic term:
\begin{align*}
\nu_2 \left( \NTerm(4n) \right) - 1
= \nu_2 \left( \frac{\HW(M(4n))}{u(4n)} - t(4n)^2 \right) - 1
= \omega(n) .
\end{align*}
\end{proof}
We see that the number of distinct prime divisors of a positive integer $n$ can be determined by a fixed-length elementary closed-form expression. Most interesting is that this expression does not reveal any information about the specific divisors of $n$. Recall that $M(n)$ \cref{TermM} is a fixed-length elementary closed-form expression and so are $\nu_2(n)$ and $\HW(n)$ as shown in \cref{SectNumberTheoreticTerms}.

The full arithmetic term for $M(n)$ is lengthy, but we include it here for completeness:
\begin{align} \label{TermMIntegerRepresentation}
M(n)
&= ( 2^{2 t(n) u(n)^2 + u(n)}-2^{u(n)} ) \cdot (2^{u(n)}+1)^{-1} + \notag \\
&\quad- (2^{2 t(n) u(n)}-1) (n 2^{2 t(n) u(n) + u(n) + 1} - n 2^{2 t(n) u(n)+1}) \ell_1
\cdot (2^{2 u(n)}-1)^{-1} (2^{2 t(n) u(n)}-1)^{-2} + \notag \\
&\quad+ (2^{3 u(n) + 1}-2^{2 u(n)+1}) (2^{2 t(n)^2 u(n)}-1) \ell_2
\cdot (2^{2 u(n)}-1)^{-3} (2^{2 t(n) u(n)}-1)^{-1} + \notag \\
&\quad+ n 2^{2 t(n) u(n)+2 u(n)+1} (2^{u(n)}-1) \ell_1 \ell_2
\cdot
    (2^{2 u(n)}-1)^{-3} (2^{2 t(n) u(n)}-1)^{-2} + \notag \\
&\quad- n^2 (2^{2 t(n) u(n)}-1) (2^{2 t(n) u(n)+u(n)}-2^{2 t(n) u(n)}) \ell_3
\cdot
    (2^{2 u(n)}-1)^{-1} (2^{2 t(n) u(n)}-1)^{-3}
+ \notag \\
&\quad- (2^{3 u(n)}-2^{2 u(n)}) (2^{2 t(n)^2 u(n)}-1) \ell_4
\cdot
(2^{2 u(n)}-1)^{-5} (2^{2 t(n) u(n)}-1)^{-1},
\end{align}
where $t(n)=n+1$, $u(n)=n+5$,
\begin{align*}
\ell_1 &= n 2^{2 t(n)^2 u(n)} - t(n) 2^{2 n t(n) u(n)} + 1, \\
\ell_2 &= n^2 2^{2 u(n)(n+2)} - (2 n^2+2n-1) 2^{2 t(n) u(n)} + t(n)^2 2^{2 n u(n)}-2^{2 u(n)}-1 , \\
\ell_3 &= n^2 2^{2 t(n) u(n)(n+2)} - (2 n^2+2n-1) 2^{2 t(n)^2 u(n)} + t(n)^2 2^{2 n t(n) u(n)} - 2^{2 t(n) u(n)}-1 , \\
\ell_4 &= (6 n^4+12 n^3-6 n^2-12n+11) 2^{2 u(n)(n+2)}
+ (-4 n^4-12 n^3-6 n^2+12n+11) 2^{2 t(n) u(n)} + \\
&\quad + (-4 n^4-4 n^3+6 n^2-4 n+1) 2^{2 u(n)(n+3)}
+ t(n)^4 2^{2 n u(n)} -2^{6 u(n)} - 11 \cdot 2^{4 u(n)} - 11 \cdot 2^{2 u(n)}-1 + \\
&\quad + n^4 2^{2 u(n)(n+4)} .
\end{align*}

\begin{remark}
The above arithmetic term for $M(n)$ can be generated by the Maple source code included in Appendix \cref{appendix:maple} and the SageMath source code included in Appendix \cref{appendix:sagesourcecodesarithmeticterms}.
\end{remark}

\section{The prime-counting function}
The \textbf{prime-counting function}, $\pi(n)$, is defined for natural numbers $n$ and returns the number of primes less than or equal to $n$.

\begin{lemma} \label{proof:piomega}
$ \forall n \in \mathbb{N}, \quad \omega(n!) = \pi(n) .$
\end{lemma}
\begin{proof}
The proof is trivial, though we note that the factorial function $n! = n (n-1)!$ defines $0! = 1$.
\end{proof}

\begin{theorem} \label{proof:primepi}
For all $n \in \mathbb{N}$, the number of primes less than or equal to $n$, $\pi(n)$, is given by the arithmetic term:
\begin{align} \label{TermPi}
\pi(n) = \nu_2\left( \frac{\HW(M(4(n!)))}{u(4(n!))} - t(4(n!))^2 \right) - 1 ,
\end{align}
where $t(n)=n+1$, $u(n)=n+5$, and $M(n)$ is the arithmetic term defined in \cref{TermM}.
\end{theorem}
\begin{proof}
Let $n \in \mathbb{N}$. By \cref{proof:piomega}, $\omega(n!) = \pi(n)$ which is obvious. Applying \cref{proof:omega}, we obtain:
\begin{align*}
\omega(n!) = \nu_2\left( \frac{\HW(M(4(n!)))}{u(4(n!))} - t(4(n!))^2 \right) - 1 = \pi(n) .
\end{align*}
\end{proof}

\section{An exponential Diophantine equation for the n-th prime} \label{section:nthprimediophantine}
Our final task is to find an arithmetic term for the $n$-th prime number, $p(n)$. We will commence by defining a system of exponential Diophantine equations derived from the results in previous sections.

Of considerable importance will be the arithmetic term $M(n)$ \cref{TermM}, which we used in previous sections to develop our arithmetic terms for $\NTerm(n)$, $\omega(n)$, and $\pi(n)$. We will require the single-fold relation $m = M(4 n)$, which will be written as:
\begin{align*}
    E_{M4}(n, [0], m) = 0 .
\end{align*}
To construct a square which defines this relation, we must find a \textbf{normalized rational form} for $M(n)$, which is an expression $M(n) = \frac{\Numer(n)}{\Denom(n)}$ such that $\Numer(n)$ and $\Denom(n)$ are arithmetic terms. Normalizing $M(n)$ involves finding a common denominator for all terms in $M(n)$ and an amplifier to ensure that both $\Numer(n)$ and $\Denom(n)$ do not contain any terms with rational coefficients.
\begin{remark}
A normalized rational form $M(n) = \frac{\Numer(n)}{\Denom(n)}$ can be calculated using the Maple source code included in Appendix \cref{appendix:maple} via Maple's built-in $\operatorname{normal}(\cdot)$ function.
\end{remark}

\begin{lemma} \label{proof:msingle-fold}
\begin{align*}
\forall (x_1,x_2) \in \mathbb{N}^2, \;
x_2 = M(4 x_1) \iff E_{M4}(x_1, [0], x_2) = 0
\iff (x_2 \cdot \Denom(4 x_1) - \Numer(4 x_1))^2 = 0 ,
\end{align*}
where $M(n)$ is the arithmetic term defined in \cref{TermM} and $\Numer(n)$, $\Denom(n)$ are arithmetic terms such that $M(n) = \frac{\Numer(n)}{\Denom(n)}$ for all $n \in \mathbb{N}$.
\end{lemma}

\begin{lemma} \label{proof:primeinequality}
$ \forall n \in \Z^+, \quad
p(n) = | \{ a \in \{ 0,\ldots,n^2 \} : \NTerm(4(a!)) \leq 2^n \} | . $
\end{lemma}
\begin{proof}
Let $n \in \Z^+$. It follows from \cref{proof:primepi} that
\begin{align*}
\forall a \in \mathbb{N}, \quad
\NTerm(4(a!)) = 2^{\omega(a!) + 1} = 2^{\pi(a) + 1} .
\end{align*}
As stated by Jones in \cite{jones1975formula}, $p(n)$ is equal to the number of solutions $a \in \{ 0, \ldots, n^2 \}$ to the inequality
\begin{align*}
\pi(a) < n .
\end{align*}
We can avoid direct comparisons with $\pi(a)$ by observation that the sets
\begin{align*}
\{ a \in \{ 0,\ldots,n^2 \} : 2^{\pi(a) + 1} \leq 2^n \}
\end{align*}
and
\begin{align*}
\{ a \in \{ 0,\ldots,n^2 \} : \pi(a) < n \}
\end{align*}
coincide. The lemma follows immediately.
\end{proof}

\begin{lemma} \label{proof:nthprimesystem}
Let $\vec{x} = (x_1,x_2,x_3,x_4)$ and let $\vec{y} = (y_1, y_2, \ldots, y_{25})$. Then for all $n \in \Z^+$, the number of solutions $(\vec{x},\vec{y}) \in \mathbb{N}^{29}$ to the system of exponential Diophantine equations
\begin{align} \label{eq:nthprimesystem}
E_{!}(x_1, [13], x_2) &= 0 , \notag \\
E_{M4}(x_2, [0], x_3) &= 0 , \notag \\
E_{\HW}(x_3, [12], x_4) &= 0 , \notag \\
\frac{x_4}{4 x_2 + 5} - (4 x_2 + 1)^2 &\leq 2^n ,
\end{align}
equals the $n$-th prime number, $p(n)$.
\end{lemma}
\begin{proof}
The first equation, $E_{!}(x_1, [13], x_2) = 0$, defines $x_2 = x_1 !$ single-fold.

The second equation, $E_{M4}(x_2, [0], x_3) = 0$, defines $x_3 = M(4 x_2) = M(4(x_1 !))$ single-fold.

The third equation, $E_{\HW}(x_3, [12], x_4) = 0$, defines $x_4 = \HW(x_3) = \HW(M(4(x_1 !)))$ single-fold.

From \cref{proof:ntermarithmetic}, we see that our final inequality, $\frac{x_4}{4 x_2 + 5} - (4 x_2 + 1)^2 \leq 2^n$, can be rewritten as
\begin{align*}
\frac{x_4}{4 x_2 + 5} - (4 x_2 + 1)^2
= \NTerm(x_3) = \NTerm(4(x_1 !)) \leq 2^n ,
\end{align*}
and by \cref{proof:primeinequality}, we have
\begin{align*}
\forall n \in \Z^+, \quad
p(n) = | \{ x_1 \in \mathbb{N} : \NTerm(4(x_1 !)) \leq 2^n \} | .
\end{align*}
Thus, the number of solutions to the system is $p(n)$.

Finally, the number of quantified variables in the system is
\begin{align*}
[x_1,x_2,x_3,x_4] + [13] + [12] = [4] + [25] = [29] .
\end{align*}
Therefore, for all $n \in \Z^+$, the number of solutions to the system is generated by the set $\{ (\vec{x},\vec{y}) \in \mathbb{N}^{29} \}$.
\end{proof}

\subsection{Building the equation}
Now that we have a system of exponential Diophantine equations with the property that the number of solutions for all $n \in \Z^+$ is the $n$-th prime number $p(n)$ \cref{eq:nthprimesystem}, our next objective is to construct a single exponential Diophantine equation $F(n,\vec{x}) = 0$ whose solution set is exactly the solution set of the system. For $n \in \Z^+$, the number of solutions $\vec{x} \in \mathbb{N}^k$ will be exactly $p(n)$.

During the construction of the equation $F(n,\vec{x}) = 0$, we will relabel some of the $x_i$ variables in the original system to make it easier to distinguish which variables arise from certain single-fold relations.
\begin{remark}
The benefits of this relabeling will be more clear in \cref{section:monomials}, where we will reduce the total number of monomials in the equation.
\end{remark}
We commence now constructing the exponential Diophantine equation $F=0$. First, we introduce
\begin{align*}
F(a, f_1, [13]) = E_{!}(a, [13], f_1) = 0 ,
\end{align*}
which enforces the relation $f_1 = a!$ in our equation.

Second, we introduce a new variable $m$, so that $m = M(4 f_1)$. This is accomplished by adding $E_{M4}$, which is a sum of squares. Our equation is now
\begin{align*}
F(a, f_1, m, [13]) = E_{!}(a, [13], f_1) + E_{M4}(f_1, [0], m) = 0 .
\end{align*}
The expansion of $E_{M4}(f_1, [0], m)$ contains monomials with non-linear exponents. In particular, several monomials contain the terms $f_1^2$ and $f_1^3$ in their exponents. This conflicts with the hypercube method requirement that all variables in exponents be simply occurring, hence linear functions in the displayed unknowns. To workaround this, we must introduce new variables $f_2 = f_1^2$, $f_3 = f_1^3$ and update our equation to be
\begin{align*}
F(a, f_1, f_2, f_3, m, [13]) = E_{!}(a, [13], f_1) + (f_2 - f_1^2)^2 + (f_3 - f_1 f_2)^2 + E_{M4}(f_1, f_2, f_3, [0], m) = 0 ,
\end{align*}
where $E_{M4}(f_1, f_2, f_3, [0], m)$ is identical to $E_{M4}(f_1, [0], m)$, except that all instances of $f_1^2$ and $f_1^3$ in the exponents have been replaced by $f_2$ and $f_3$ respectively.

Third, we introduce a new variable $b$, so that $b = \HW(m)$. This is accomplished by adding the term $E_{\HW}(m, [12], b)$, which is again a sum of squares. Now, we have
\begin{align*}
F(a, f_1, f_2, f_3, m, b, [25]) &= E_{!}(a, [13], f_1) + (f_2 - f_1^2)^2 + (f_3 - f_1 f_2)^2 + E_{M4}(f_1, f_2,  f_3, [0], m) \\ &\quad + E_{\HW}(m, [12], b) = 0 .
\end{align*}

It remains only to enforce the inequality $\NTerm(4a!) \leq 2^n$. In our equation, we have that $\NTerm(4(a!)) = \frac{b}{4 f_1 + 5} - (4 f_1 + 1)^2$. We introduce now a new variable $d$, so that
\begin{align*}
F(n, a, f_1, f_2, f_3, m, b, d, [25])
&= E_{!}(a, [13], f_1) + (f_2 - f_1^2)^2 + (f_3 - f_1 f_2)^2 + E_{M4}(f_1, f_2, f_3, [0], m) \\
&\quad + E_{\HW}(m, [12], b) + \left(\frac{b}{4 f_1 + 5} - (4 f_1 + 1)^2 + d - 2^n \right)^2 = 0\end{align*}
To remove the denominator in the final square term, we multiply the inner terms by $(4 f_1 + 5)$. This leads to
\begin{align} \label{eq:F}
F(n, a, f_1, f_2, f_3, m, b, d, [25])
&= E_{!}(a, [13], f_1) + (f_2 - f_1^2)^2 + (f_3 - f_1 f_2)^2 + E_{M4}(f_1, f_2, f_3, [0], m) \notag \\
&\quad + E_{\HW}(m, [12], b) + \left(b + (4 f_1 + 5) \left( -(4 f_1 + 1)^2 + d - 2^n\right) \right)^2 = 0 .
\end{align}

Finally, we relabel the variables $(a, f_1, f_2, f_3, m, b, d)$ over the tuple $\vec{x} = (x_1,x_2,\ldots,x_{7})$ and write
\begin{align} \label{eq:Fxy}
    F(n, \vec{x}, [25]) = 0 \iff F(n, a, f_1, f_2, f_3, m, b, d, [25]) = 0 .
\end{align}

\begin{theorem} \label{proof:primeeq1}
Let $\vec{x} = (x_1,x_2,\ldots,x_7)$ and let $\vec{y} = (y_1,y_2,\ldots,y_{25})$. For all $n \in \Z^+$, the $n$-th prime $p(n)$, is given by the number of solutions $(\vec{x},\vec{y}) \in \mathbb{N}^{32}$ to the exponential Diophantine equation $F(n,\vec{x},[25]) = 0$.
\end{theorem}
\begin{proof}
The equation $F(n, \vec{x}, [25])=0$ translates the system of exponential Diophantine equations proved in \cref{proof:nthprimesystem} to a single exponential Diophantine equation with relabeling of certain variables. The total number of quantified variables in the equation is
\begin{align*}
[x_1,x_2,x_3,x_4,x_5,x_6,x_7] + [25] = [7] + [25] = [32] .
\end{align*}
The final equation $F(n, \vec{x}, [25])=0$ is a sum of squares of the equations defined in \cref{proof:nthprimesystem}, which are all single-fold relations. Furthermore, the values of all single-fold relations in the set of equations are uniquely determined by the variable $x_1$. Hence, by \cref{proof:sumofsquares}, the solution set of the system coincides with the solution set $\{ (\vec{x},\vec{y}) = (x_1, x_2, \ldots, x_7, y_1, y_2, \ldots, y_{25}) \in \mathbb{N}^{32} \} : F(n,\vec{x},[25]) = 0$.
\end{proof}


\section{Reducing the number of monomials} \label{section:monomials}
Our exponential Diophantine equation $F(n, \vec{x}, [25]) = 0$ developed in \cref{section:nthprimediophantine} could theoretically be used to construct an arithmetic term for the $n$-th prime.

Since the hypercube method is typically defined to accept a single tuple of variables $\vec{x}$, we will rewrite $F(n,\vec{x},[25]) = 0$ as $F(n,\vec{x}) = 0$ by extending the tuple $\vec{x} = (x_1,x_2,\ldots,x_{32})$ and then relabeling the variables $\vec{y} = (y_1,y_2,\ldots,y_{25})$ over the appended $x_i$ variables $(x_{8},x_{9},\ldots,x_{32})$ respectively. Thus,
\begin{align*}
     & F(n,\vec{x}) = F(n,x_1,x_2,\ldots,x_{32}) = 0 \\
     & \iff
     F(n,x_1,x_2,\ldots,x_7,[25]) = F(n,x_1,x_2,\ldots,x_7,y_1,y_2,\ldots,y_{25}) = 0 .
\end{align*}
A problem is that we have not yet calculated the bounds for the quantified variables $\vec{x}$ as a function of $n \in \mathbb{Z}^+$. Without these bounds, we cannot yet determine appropriate arithmetic terms $t(n)$ and $u(n)$, which are the parameters for hypercube method and whose growth determines if the resulting arithmetic term is universally applicable for all $n \in \Z^+$. Nonetheless, as it is theoretically feasible, we will state a partial result:
\begin{theorem}
For all $n \in \Z^+$, the $n$-th prime number, $p(n)$, is given by an arithmetic term:
\begin{align*}
p(n) = \frac{\HW(Q(n))}{u(n)} - t(n)^{32} ,
\end{align*}
where $t(n), u(n)$ are arithmetic terms chosen in accordance with the hypercube method and $Q(n)$ is the arithmetic term constructed by the hypercube method from the monomial expansion of the exponential Diophantine equation $F(n,\vec{x})=0$ with the parameters $k$, $t(n)$, and $u(n)$.

Let $m_i(n,\vec{x})$ represent a monomial in the monomial expansion of $F(n,\vec{x}) = \sum_{i=0}^j m_i(n,\vec{x}) = 0$, where the monomial ordering scheme is arbitrary and the number of monomials is $j=10102$.

The arithmetic term $Q(n)$ arises as the sum of a single $\CTerm_{32}(m_i(n,\vec{x}),t(n),u(n))$ term and $10099$ distinct $\ATerm_{32}(m_i(n,\vec{x}),t(n),u(n))$ terms. Every $\ATerm_{32}(m_i(n,\vec{x}),t(n),u(n))$ term is the product of $32$ $G_r(b,t(n))$ terms, where $r \in \{0,\ldots,31\}$ and $b \in \mathbb{N}$ is determined by $n$ and $\vec{x}$.
\end{theorem}
Unfortunately, this arithmetic term for $p(n)$ is far too lengthy: Merely writing its $\CTerm_k(\cdot)$ term and many $\ATerm_k(\cdot)$ terms would occupy roughly $50$ pages in this paper. Consequently, establishing the bounds for the arithmetic terms $t(n)$, $u(n)$ would be an extremely tedious exercise and also cumbersome to verify. We opt instead to find an arithmetic term for $p(n)$ that can be written down completely on a few pages as a sum of its $\CTerm_k(\cdot)$ term and $\ATerm_k(\cdot)$ terms. This necessitates reducing the number of monomials in our exponential Diophantine equation $F(n,\vec{x},[25])=0$ from $10102$ down to fewer than $500$. Such term will allow our results to be more easily analyzed and shared.

It is possible to significantly reduce the number of monomials by introducing new variables and breaking the large square relations down into smaller squares. However, each additional variable introduced increases the magnitude of the resulting arithmetic term, creating a trade-off. While the final $p(n)$ term we present in \cref{section:nthprimeterm} requires fewer $\CTerm_k(\cdot)$ and $\ATerm_k(\cdot)$ terms to write down, the $32$ variable version is in actuality much smaller from a computational perspective.

\begin{remark}
 To clarify: The number of monomials roughly corresponds to the number of $\ATerm_k(\cdot)$ terms comprising the sum of the final arithmetic term, while the number of variables corresponds to the number of $G_r(\cdot)$ terms that each $\ATerm_k(\cdot)$ term is the product of \cref{TermA}. Thus, for each new variable introduced, every $\ATerm_k(\cdot)$ term will be multiplied by an additional $G_r(\cdot)$ term. These $G_r(\cdot)$ multiplications far outweigh the summing of additional $\ATerm_k(\cdot)$ terms arithmetically.
\end{remark}

 For the monomial reduction, we will require many additional variables. To assist the reader, we begin again by rewriting $F(n, \vec{x}, [25]) = 0$ \cref{eq:F} without the initial $\vec{x}$ variable relabelings, so that:
\begin{align*}
F(n, a, f_1, f_2, f_3, m, b, d, [25]) = 0
\iff F(n, \vec{x}, [25]) = F(n, x_1, x_2, \ldots, x_7, [25]) = 0 .
\end{align*}
In this equation, the vast majority of monomials arise from the expansion of the square term $E_{M4}$, which produces a staggering $9920$ of the $10102$ total monomials. Our goal now is to reduce the number of monomials in the expansion of this square by splitting it up into a sum of smaller squares.

First, we define some temporary \textit{non-quantified} variables:
\begin{align*}
& u = 4 f_1 + 5 ,
\quad
t = 4 f_1 + 1 ,
\quad
t_1 = t - 1,
\quad
q_1 = 2^{2u} ,
\quad
q_2 = 2^{2tu} .
\end{align*}

Now, we introduce the \textit{quantified} variables
\begin{align*}
f_4 = 4 f_1 + 1, \quad
q_{1,1} = q_1 - 1 , \quad
q_{2,1} = q_2 - 1 ,
\end{align*}
as the sum of squares
\begin{align*}
(f_4 - 4 f_1 - 1)^2
+ (q_{1,1} - q_1 + 1)^2
+ (q_{2,1} - q_2 + 1)^2 = 0 .
\end{align*}
Importantly, the variable $f_4$ will be quantified outside the expression and passed in, so that it can be used elsewhere.

Next, we define additional \textit{quantified} variables for the $G_r(\cdot)$ terms \cref{TermG} comprising the product expansion of our $\ATerm_k(\cdot)$ terms \cref{TermA}:
\begin{align*}
\begin{array}{lll}
g_{0,1} = G_0(q_1, t) , &
g_{0,2} = G_0(q_2, t) , &
g_{1,2} = G_1(q_2, t) , \\
g_{2,2} = G_2(q_2, t) , &
g_{2,1} = G_2(q_1, t) , &
g_{4,1} = G_4(q_1, t) .
\end{array}
\end{align*}
Expanding these, we have
\begin{align*}
g_{0,1} &= (q_1^t - 1) / q_{1,1} = q_{2,1} / q_{1,1}, \\
g_{0,2} &= (q_2^t - 1) / q_{2,1}, \\
g_{1,2} &= q_2 \left( t_1 q_2^t - f_4 q_2^{t_1} \right) / q_{2,1}^2 , \\
g_{2,1} &= q_1 \left(t_1^2 q_1^{t_1 + 2} - (2 t_1^2 + 2 t_1 - 1) q_1^t + f_4^2 q_1^{t_1} - q_1 - 1
\right) / q_{1,1}^3, \\
g_{2,2} &= q_2 \left( t_1^2 q_2^{t_1 + 2} - (2 t_1^2 + 2 t_1 - 1) q_2^t + f_4^2 q_2^{t_1} - q_2 - 1
\right) / q_{2,1}^3, \\
g_{4,1} &=
q_1^{t_1 + 2} (-4 \cdot t_1^{4} - 12 t_1^{3} - 6 t_1^{2} + 12 t_1 + 11) / q_{1,1}^{5} \\
&\quad + q_1^{t_1 + 3} (6 t_1^{4} + 12 t_1^{3} - 6 t_1^{2} - 12 t_1 + 11) / q_{1,1}^{5} \\
&\quad + q_1^{t_1 + 4} (-4 \cdot t_1^{4} - 4 t_1^{3} + 6 t_1^{2} - 4 t_1 + 1) / q_{1,1}^{5} \\
&\quad + q_1 (f_4^{4} q_1^{t_1} + t_1^{4} q_1^{t_1} - q_1^{3} - 11 q_1^{2} - 11 q_1 + 3) / q_{1,1}^{5} .
\end{align*}

Putting these together, the sum of squares defining our $G_r(\cdot)$ terms is
\begin{align*}
& \left( g_{0,1} q_{1,1} - (q_1^t - 1) \right)^2 \\
&+ \left( g_{0,2} q_{2,1} - (q_2^t - 1) \right)^2 \\
&+ \left( g_{1,2} q_{2,1}^2 - q_2 \left( t_1 q_2^t - f_4 q_2^{t_1} \right) \right)^2 \\
&+ \left( g_{2,1} q_{1,1}^3 - q_1 \left(t_1^2 q_1^{t_1 + 2} - (2 t_1^2 + 2 t_1 - 1) q_1^t + f_4^2 q_1^{t_1} - q_1 - 1
\right) \right)^2 \\
& + \left(g_{2,2} q_{2,1}^3 - q_2 \left( t_1^2 q_2^{t_1 + 2} - (2 t_1^2 + 2 t_1 - 1) q_2^t + f_4^2 q_2^{t_1} - q_2 - 1
\right) \right)^2 \\
&+ \left( g_{4,1} q_{1,1}^{5} - \numerator(G_4(q_1, t)) \right)^2 = 0 ,
\end{align*}
where the $G_r(\cdot)$ functions are replaced by the relations defined above for the $g_{i,j}$.

Next, we introduce a \textit{non-quantified} variable for the common factor used in the calculations of $\ATerm_k(\cdot)$ as
\begin{align*}
h = 1 - 2^u = -2^u + 1 .
\end{align*}

For additional clarity, we record the following additional \textit{non-quantified} variables, which each correspond to a specific $\ATerm_k(\cdot)$ term in the sum $M(\cdot)$:
\begin{align*}
A_{4,0} &= h g_{4,1} g_{0,2} , \quad
A_{2,0} = h g_{2,1} g_{0,2} , \quad
A_{2,1} = h g_{2,1} g_{1,2} , \quad
A_{0,2} &= h g_{0,1} g_{2,2} , \quad
A_{0,1} = h g_{0,1} g_{2,1} .
\end{align*}

Finally, we introduce a \textit{quantified} variable for our $\CTerm_k(\cdot)$ term \cref{TermC} as
\begin{align*}
C = \frac{2^u (2^{2ut^2} - 1)}{2^u + 1}
= \frac{2^u (q_2^t - 1)}{2^u + 1} ,
\end{align*}
which corresponds to the square
\begin{align*}
\left( C \cdot (2^u + 1) - (2^u - t + 1) (q_2^t - 1) \right)^2 = 0.
\end{align*}
Altogether, we have
\begin{align*}
M(4 f_1)
&= C + A_{4,0}
- 2 A_{2,0}
- 2 (4 f_1) A_{2,1}
+ (4 f_1)^2 A_{0,2}
+ 2 (4 f_1) A_{0,1} ,
\end{align*}
which is defined by the square
\small
\begin{align} \label{EqEM1}
E_{M4}(f_1, f_2, f_3, f_4, [9], m)
= \left(
m -
(C + A_{4,0}
- 2 A_{2,0}
- 2 (4 f_1) A_{2,1}
+ (4 f_1)^2 A_{0,2}
+ 2 (4 f_1) A_{0,1} )
\right)^2
= 0 .
\end{align}
\normalsize
The number of \textit{quantified} variables is
\begin{align*}
& [q_{1,1},q_{2,1}]
+ [g_{0,1},g_{2,1},g_{4,1},g_{0,2},g_{1,2},g_{2,2}]
+ [C] \\
&= [2] + [6] + [1] \\
&= [9] .
\end{align*}
In the lemma below, we relabel the variables $(f_1,f_2,f_3,f_4,m)$ over the tuple $(x_1,x_2,x_3,x_4,x_5)$ respectively. The remaining \textit{quantified variables} are relabeled over a tuple $\vec{y} = (y_1,y_2,\ldots,y_{9})$, where the order in which the relabeling is performed is consistent with their respective definitions:
\begin{lemma} \label{proof:msingle-fold2}
For all $(x_1, x_2, x_3, x_4, x_5) \in \mathbb{N}^5$ such that $x_2 = x_1^2$, $x_3 = x_1^3$, and $x_4 = 4 x_1 + 1$, we have that
\begin{align*}
x_5 = M(4 x_1)
\iff E_{M4}(x_1, x_2, x_3, x_4, [9], x_5) = 0 .
\end{align*}
where $M(n)$ is the arithmetic term defined in \cref{TermM} for $n \in \mathbb{N}$.
\end{lemma}
Using this new version of $E_{M4}$, whose squares produce far fewer monomials than the original relation, we define a new exponential Diophantine equation whose solution count in the natural numbers is the $n$-th prime number $p(n)$, as:
\begin{align} \label{eq:F1}
\hat{F}(n, a, f_1, f_2, f_3, f_4, m, b, d, [34])
&= E_{!}(a, [13], f_1)
+ (f_2 - f_1^2)^2
+ (f_3 - f_1 f_2)^2
+ (f_4 - 4 f_1 - 1)^2
\notag \\
&\quad
+ E_{M4}(f_1, f_2, f_3, f_4, [9], m)
+ E_{\HW}(m, [12], b)
\notag \\
&\quad
+ \left(b + (4 f_1+5) \left( -f_4^2 + d - 2^n\right) \right)^2
= 0 ,
\end{align}
where the number of quantified variables is
\begin{align*}
[a, f_1, f_2, f_3, f_4, m, b, d] + [13] + [9] + [12]
= [8] + [34] = [42] .
\end{align*}
After relabeling the variables $(a, f_1, f_2, f_3, f_4, m, b, d)$ over the tuple $\vec{x} = (x_1,x_2,\ldots,x_{8})$, one has:
\begin{align} \label{eq:F1xy}
    \hat{F}(n,\vec{x},[34]) = 0
    \iff
    \hat{F}(n, a, f_1, f_2, f_3, f_4, m, b, d, [34]) = 0 .
\end{align}
This equation contains $498$ monomials in $42$ quantified variables. For its full monomial expansion, please consult Appendix \cref{appendix:F1monomials}.

\begin{theorem} \label{proof:primeeq2}
Let $\vec{x} = (x_1,x_2,\ldots,x_7)$ and let $\vec{y} = (y_1,y_2,\ldots,y_{34})$. For all $n \in \Z^+$, the $n$-th prime number $p(n)$, is given by the number of solutions $(\vec{x},\vec{y}) \in \mathbb{N}^{42}$ to the exponential Diophantine equation $\hat{F}(n,\vec{x},[34]) = 0$.
\end{theorem}
\begin{proof}
The equation
\begin{align*}
    \hat{F}(n,\vec{x},[34]) = \hat{F}(n,x_1,x_2,\ldots,x_8,[34]) = 0
\end{align*}
has the same number of solutions as the equation proved in \cref{proof:primeeq1}, which is:
\begin{align*}
    F(n,x_1,x_2,\ldots,x_7,[25]) = 0 .
\end{align*}
\end{proof}

\section{Single-fold bounds} \label{section:single-foldbounds}
To write down an explicit arithmetic term for $p(n)$, we must find suitable arithmetic terms $t(n),u(n)$ that allow us to apply the hypercube method to count the number of solutions to our exponential Diophantine equation $\hat{F}(n,\vec{x})=0$. The first step in this process is to establish the bounds for all quantified variables in our single-fold relations.

\begin{lemma} \label{proof:divsingle-foldbounds}
The relation
\begin{align*}
x_3 = \floor{x_1 / x_2} \iff
E_{/}(x_1, x_2, [2], x_3) = 0 \iff \exists \vec{y} \in \mathbb{N}^2 : (x_1 - x_2 x_3 - y_1)^2 + (y_1 + y_2 + 1 - x_2)^2 = 0 ,
\end{align*}
has the bounds:
\begin{align*}
\begin{array}{lll}
y_1 < x_2 ,
& y_2 < x_2 ,
& x_3 \leq x_1 .
\end{array}
\end{align*}
\end{lemma}

\begin{lemma} \label{proof:modsingle-foldbounds}
The relation
\begin{align*}
x_3 &= x_1 \bmod x_2
\iff E_{\bmod}(x_1, x_2, [2], x_3) = 0 \\
&\iff \exists \vec{y} \in \mathbb{N}^2 : (x_1 - x_2 y_1 - x_3)^2 + (x_3 + y_2 + 1 - x_2)^2 = 0 ,
\end{align*}
has the bounds:
\begin{align*}
\begin{array}{lll}
x_3 < x_2 ,
& y_2 < x_2 ,
& y_1 < x_1 .
\end{array}
\end{align*}
\end{lemma}

\begin{lemma} \label{proof:midsingle-foldbounds}
The relation
\begin{align*}
& x_2 \mid x_1
\iff E_{|}(x_1, [1], x_2) = 0
\iff \exists y_1 \in \mathbb{N} : (x_1 - x_2 y_1)^2 = 0 ,
\end{align*}
has the bound:
\begin{align*}
   y_1 \leq x_1 .
\end{align*}
\end{lemma}

\begin{lemma} \label{proof:nmidsingle-foldbounds}
The relation
\begin{align*}
& x_2 \nmid x_1
\iff E_{\nmid}(x_1, [3], x_2) = 0
\iff \exists \vec{y} \in \mathbb{N}^3 : (x_1 - x_2 y_1 - y_2 - 1)^2 + (y_2 + y_3 + 2 - x_2)^2 = 0 ,
\end{align*}
has the bounds:
\begin{align*}
\begin{array}{lll}
y_2 < x_2 ,
& y_3 < x_2 ,
& y_1 < x_1
\end{array}
\end{align*}
\end{lemma}

\begin{lemma} \label{proof:nusingle-foldbounds}
The relation
\begin{align*}
x_2 = \nu_2(x_1)
\iff E_{\nu}(x_1, [4], x_2) = 0
\iff \exists \vec{y} \in \mathbb{N}^4 : E_\nmid (x_1, [3], 2^{x_2+1}) + E_| (x_1, [1], 2^{x_2}) = 0 ,
\end{align*}
has the bounds:
\begin{align*}
\begin{array}{lllll}
y_1 < x_1 ,
& y_2 <  2 x_1 ,
& y_3 < 2 x_1 ,
& y_4 < x_1 + 1 ,
& x_2 < x_1 .
\end{array}
\end{align*}
\end{lemma}
\begin{proof}
We have
\begin{align*}
E_{\nu} (x_1, y_1, y_2, y_3, y_4, x_2)
&= E_{\nmid} (x_1, y_1, y_2, y_3, 2^{x_2+1}) + E_{|}(x_1, y_4 , 2^{x_2}) \\
&= (x_1 -  2^{x_2+1} y_1 - y_2 - 1)^2 + (y_2 + y_3 + 2 -  2^{x_2+1})^2 + (x_1 - 2^{x_2} y_4)^2 .
\end{align*}
We observe that $2^{x_2} \leq x_1$ and $2^{x_2 + 1} \leq 2 x_1$. The bounds are
\begin{align*}
y_1 < x_1, \quad
y_2 < 2^{x_2 + 1} \leq 2 x_1, \quad
y_3 < 2^{x_2 + 1} \leq 2 x_1, \quad
y_4 \leq x_1 < x_1 + 1 , \quad
x_2 < x_1.
\end{align*}
so
\begin{align*}
y_1 < x_1,\quad
y_2 <  2 x_1, \quad
y_3 < 2 x_1, \quad
y_4 < x_1 + 1 , \quad
x_2 < x_1.
\end{align*}
\end{proof}

\begin{lemma} \label{proof:binomialcoeffsingle-foldbounds}
The relation
\begin{align*}
& x_3 = \binom{x_1}{x_2}
\iff E_{\binom{\#}{\#}}(x_1, x_2, [7], x_3) = 0
\end{align*}
\begin{align*}
\iff \exists \vec{y} \in \mathbb{N}^7 : \quad
& (y_1 - (2 x_1^3 + 8 x_1^2 + 2 x_1 x_2 + 12 x_1 + 4 x_2 + 8 ))^2
+ (y_2 - (2 x_1^2 + 8 x_1 + 8))^2 \\
& + E_{/} (2^{y_1}, 2^{y_2} - 2^{2 x_1 + 4} - 1, [2], y_3)
+ E_{\bmod} (y_3, 2^{2 x_1 + 4}, [2] , x_3) = 0 ,
\end{align*}
has the bounds:
\begin{align*}
\begin{array}{lll}
y_1 < 28 x_1^3 + 9,
& y_2 < 10 x_1^2 + 9,
& y_3 < 2^{ 28 x_1^3 + 9 },
\\
y_4 < 2^{10 x_1^2 + 8} ,
& y_5 < 2^{10 x_1^2 + 8} ,
& y_6 < 2^{ 28 x_1^3 + 9 } , \\
y_7 < 2^{2 x_1 + 4} ,
& x_3 < 2^{2 x_1 + 4} .
\end{array}
\end{align*}
\end{lemma}
\begin{proof}
We start with
\begin{align*}
y_1  = 2 x_1^3 + 8 x_1^2 + 2 x_1 x_2 + 12 x_1 + 4 x_2 + 8
\end{align*}

We know that $x_2 \leq x_1$, so we can majorize $x_2$ with $x_1$. Hence
\begin{align*}
y_1 &\leq 2 x_1^3 + 8 x_1^2 + 2 x_1^2 + 12 x_1 + 4 x_1 + 8 , \\
y_1 &\leq 2x_1^3 + 10 x_1^2 + 16 x_1 + 8 , \\
y_1 &\leq  28 x_1^3 + 8 < 28 x_1^3 + 9 .
\end{align*}

As a principle, we will always let a constant term inside
to be sure that the strict inequality works also for $x_1 = 0$.

We also have
\begin{align*}
y_2 = 2 x_1^2 + 8 x_1 + 8 \leq 10 x_1^2 + 8 < 10 x_1^2 + 9 .
\end{align*}

By \cref{proof:divsingle-foldbounds}, the relation
\begin{align*}
E_{/} (2^{y_1}, 2^{y_2} - 2^{2 x_1 + 4} - 1, y_4, y_5, y_3) = 0 ,
\end{align*}
has the bounds:
\begin{align*}
y_4 &< 2^{y_2} - 2^{2 x_1 + 4} - 1 < 2^{10 x_1^2 + 8} , \\
y_5 &< 2^{y_2} - 2^{2 x_1 + 4} - 1 < 2^{10 x_1^2 + 8} , \\
y_3 &\leq  2^{y_1} < 2^{ 28 x_1^3 + 9 } .
\end{align*}

Finally, by \cref{proof:modsingle-foldbounds}, the relation
\begin{align*}
E_{\bmod} (y_3, 2^{2 x_1 + 4}, y_5, y_6 , x_3) = 0 ,
\end{align*}
has the bounds:
\begin{align*}
y_6 < 2^{ 28 x_1^3 + 9 } , \quad
y_7 < 2^{2 x_1 + 4} , \quad
x_3 < 2^{2 x_1 + 4} .
\end{align*}
\end{proof}

\begin{lemma} \label{proof:factorialsingle-foldbounds}
The relation
\begin{align*}
x_2 = x_1! \iff E_{!}(x_1, [13], x_2) = 0
\end{align*}
\begin{align*}
\iff \exists \vec{y} \in \mathbb{N}^{13} : \;
& (y_1 - x_1^2)^2
+ (y_2 - 2^{3 y_1})^2
+ (y_3 - x_1 y_1)^2
+ E_{\binom{\#}{\#}}(y_2, x_1, [7], y_4)  \\
&+ E_{/}(2^{3 y_3}, y_4, [2], x_2)
= 0 ,
\end{align*}
has the bounds:
\begin{align*}
\begin{array}{lll}
y_1 < x_1^2 + 1 ,
& y_2 < 2^{3 y_1 + 1} ,
& y_3 < x_1^3 + 1 , \\
y_4 < 2^{2 y_2 + 4},
& y_5 < 28 y_2^3 + 9,
& y_6 < 10 y_2^2 + 9, \\
y_7 < 2^{ 28 y_2^3 + 9 },
& y_8 < 2^{10 y_2^2 + 8} ,
& y_9 < 2^{10 y_2^2 + 8} , \\
y_{10} < 2^{ 28 y_2^3 + 9 } ,
& y_{11} < 2^{2 y_2 + 4} ,
& y_{12} < 2^{2 y_2 + 4}  , \\
y_{13} < 2^{2 y_2 + 4}  ,
& x_2 < 2^{3 y_3 + 1} .
\end{array}
\end{align*}
\end{lemma}
\begin{proof}
We start with
\begin{align*}
y_1 = x_1^2 < x_1^3 + 1 ,
y_2 = 2^{3 y_1} < 2^{3 y_1 + 1} ,
y_3 = x_1^3 < x_1^3 + 1 .
\end{align*}

By \cref{proof:binomialcoeffsingle-foldbounds}, the relation
\begin{align*}
E_{\binom{\#}{\#}}(y_2, x_1, [7], y_4) = 0 ,
\end{align*}
has the bounds:
\begin{align*}
\begin{array}{lll}
y_5 < 28 y_2^3 + 9,
& y_6 < 10 y_2^2 + 9,
& y_7 < 2^{ 28 y_2^3 + 9 }, \\
y_8 < 2^{10 y_2^2 + 8} ,
& y_9 < 2^{10 y_2^2 + 8} ,
& y_{10} < 2^{ 28 y_2^3 + 9 } , \\
y_{11} < 2^{2 y_2 + 4} ,
& y_4 < 2^{2 y_2 + 4} .
\end{array}
\end{align*}

By \cref{proof:divsingle-foldbounds}, the relation
\begin{align*}
E_{/}(2^{3 y_3}, y_4, [2], x_2) = 0 ,
\end{align*}
has the bounds:
\begin{align*}
y_{12} &< y_4 < 2^{2 y_2 + 4} , \\
y_{13} &< y_4 < 2^{2 y_2 + 4} , \\
x_2 &\leq 2^{3 y_3} < 2^{3 y_3 + 1} .
\end{align*}
\end{proof}

\begin{lemma} \label{proof:hwsingle-foldbounds}
The relation
\begin{align*}
x_2 = \HW(x_1)
\iff E_{\HW}(x_1, [12], x_2) = 0
\iff \exists \vec{y} \in \mathbb{N}^{12} :
E_{\binom{\#}{\#}}(2 x_1, x_1, [7], y_1) + E_\nu(y_1, [4], x_2) = 0 ,
\end{align*}
has the bounds:
\begin{align*}
\begin{array}{lll}
y_1 < 2^{4 x_1 + 4} ,
& y_2 < 224 x_1^3 + 9 ,
& y_3 < 40 x_1^2 + 9 , \\
y_4 < 2^{224 x_1^3 + 9} ,
& y_5 < 2^{40 x_1^2 + 8} ,
& y_6 < 2^{40 x_1^2 + 8} , \\
y_7 < 2^{224 x_1^3 + 9} ,
& y_8 < 2^{4 x_1 + 4} ,
& y_9 < 2^{4 x_1 + 4} , \\
y_{10} < 2^{4 x_1 + 5},
& y_{11} < 2^{4 x_1 + 5},
& y_{12} < 2^{4 x_1 + 4} + 1 , \\
x_2 < x_1 + 1 .
\end{array}
\end{align*}
\end{lemma}
\begin{proof}
By \cref{proof:binomialcoeffsingle-foldbounds}, the relation
\begin{align*}
E_{\binom{\#}{\#}}(2 x_1, x_1, [7], y_1) = 0,
\end{align*}
has the bounds:
\begin{align*}
    y_2 &< 28 (2 x_1)^3 + 9 < 224 x_1^3 + 9, \\
    y_3 &< 10 (2 x_1)^2 + 9 < 40 x_1^2 + 9, \\
    y_4 &< 2^{28 (2 x_1)^3 + 9} < 2^{224 x_1^3 + 9}, \\
    y_5 &< 2^{10 (2 x_1)^2 + 8} < 2^{40 x_1^2 + 8} , \\
    y_6 &< 2^{10 (2 x_1)^2 + 8} < 2^{40 x_1^2 + 8} , \\
    y_7 &< 2^{28 (2 x_1)^3 + 9} < 2^{224 x_1^3 + 9} , \\
    y_8 &< 2^{2 (2 x_1) + 4} < 2^{4 x_1 + 4} , \\
    y_1 &< 2^{2 (2 x_1) + 4} < 2^{4 x_1 + 4} .
\end{align*}
By \cref{proof:nusingle-foldbounds}, the relation
\begin{align*}
E_\nu(y_1, [4], x_2) = 0 ,
\end{align*}
has the bounds
\begin{align*}
y_9 &< y_1 < 2^{4 x_1 + 4} , \\
y_{10} &< 2 y_1 < 2^{4 x_1 + 5}, \\
y_{11} &< 2 y_1 < 2^{4 x_1 + 5}, \\
y_{12} &< y_1 + 1 < 2^{4 x_1 + 4} + 1 , \\
x_2 &= \HW(x_1) < x_1 + 1 .
\end{align*}
\end{proof}

\begin{lemma} \label{proof:msingle-foldbounds}
The relation
\begin{align*}
& x_5 = M(4 x_1)
\iff E_{M4}(x_1, x_2, x_3, x_4, [9], x_5) = 0 \\
& \iff
\left(
x_5 -
(y_9 + A_{4,0}
- 2 A_{2,0}
- 2 (4 x_1) A_{2,1}
+ (4 x_1)^2 A_{0,2}
+ 2 (4 x_1) A_{0,1} )
\right)^2
= 0 ,
\end{align*}
where
\begin{align*}
\begin{array}{lll}
x_2=x_1^2 ,
& x_3=x_1^3 ,
& x_4 = 4 x_1 + 1 ,
\end{array}
\end{align*}
$M(n)$ is the arithmetic term defined in \cref{TermM} for $n \in \mathbb{N}$, and the variables $A_{i,j}$ are defined in \cref{EqEM1}, has the bounds:
\begin{align*}
\begin{array}{ll}
y_1 < 2^{8 x_1 + 10} ,
& y_2 < 2^{32 x_2 + 48 x_1 + 10} , \\
y_3 < (4 x_1 + 1) 2^{32 x_1^2 + 48 x_1 + 10} ,
& y_4 < (4 x_1 + 1)^3 2^{32 x_1^2 + 48 x_1 + 10} , \\
y_5 < (4 x_1 + 1)^5 2^{32 x_1^2 + 48 x_1 + 10} ,
& y_6 < 2^{128 x_1^3 + 224 x_1^2 + 88 x_1 + 10}, \\
y_7 < (4 x_1 + 1)^2 2^{128 x_1^3 + 224 x_1^2 + 88 x_1 + 10},
& y_8 < (4 x_1 + 1)^3 2^{128 x_1^3 + 224 x_1^2 + 88 x_1 + 10}, \\
y_9 < 2^{128 x_1^3 + 224 x_1^2 + 92 x_1 + 15} ,
& x_5 < 2^{512 x_1^3 + 576 x_1^2 + 216 x_1 + 27} .
\end{array}
\end{align*}
\end{lemma}
\begin{proof}
Put
\begin{align*}
& t=4x_1+1, \quad
u=4x_1+5, \\
& q_1=2^{2u}, \quad
q_2=2^{2ut}. \\
& q_{1,1} = y_{1}, \quad
q_{2,1} = y_{2} .
\end{align*}

Now, consider
\begin{align*}
C &= y_{9} , \\
g_{0,1} &= y_{3} = G_0(q_1, t) , \\
g_{2,1} &= y_{4} = G_2(q_1, t) , \\
g_{4,1} &= y_{5} = G_4(q_1, t) , \\
g_{0,2} &= y_{6} = G_0(q_2, t) , \\
g_{1,2} &= y_{7} = G_1(q_2, t) , \\
g_{2,2} &= y_{8} = G_2(q_2, t) .
\end{align*}

We apply the following principle: For $b \geq 2$ and $r \geq 0$,
\begin{align*}
G_r(q, b) = \sum_{j=0}^{b-1} q^j j^r < b q^b b^r = b^{r+1} q^b .
\end{align*}

Therefore,
\begin{align*}
y_1 &= q_{1,1} = q_1 - 1 < q_1 < 2^{8 x_1 + 10} , \\
y_2 &= q_{2,1} = q_2 - 1 < q_2 < 2^{32 x_2 + 48 x_1 + 10} , \\
y_3 &= g_{0,1} = G_0(q_1,t) < t q_1^t = t q_2 < (4 x_1 + 1) 2^{2 (4 x_1 + 5) (4 x_1 + 1)} = (4 x_1 + 1) 2^{32 x_1^2 + 48 x_1 + 10} , \\
y_4 &= g_{2,1} = G_2(q_1,t) < t^3 q_1^t < (4 x_1 + 1)^3 2^{32 x_1^2 + 48 x_1 + 10} , \\
y_5 &= g_{4, 1} = G_4(q_1, t) < t^5 q_1^t < (4 x_1 + 1)^5 2^{32 x_1^2 + 48 x_1 + 10} , \\
y_6 &= g_{0,2} = G_0(q_2, t) < q_2^t <  2^{128 x_1^3 + 224 x_1^2 + 88 x_1 + 10}, \\
y_7 &= g_{1,2} = G_1(q_2,t) < t^2 q_2^t < (4 x_1 + 1)^2 2^{128 x_1^3 + 224 x_1^2 + 88 x_1 + 10}, \\
y_8 &= g_{2,2} = G_2(q_2, t) < t^3 q_2^t < (4 x_1 + 1)^3 2^{128 x_1^3 + 224 x_1^2 + 88 x_1 + 10}, \\
y_9 &= C = \frac{2^u (q_2^t - 1)}{2^u + 1} < 2^u q_2^t < 2^{128 x_1^3 + 224 x_1^2 + 92 x_1 + 15} .
\end{align*}

Next, we consider
\begin{align*}
x_5 = M(4 f_1)
&= C
+ A_{4,0}
- 2 A_{2,0}
- 2 (4 f_1) A_{2,1}
+ (4 f_1)^2 A_{0,2}
+ 2 (4 f_1) A_{0,1} \\
&< C - 2 A_{2,0} - 8 f_1 A_{2,1} .
\end{align*}

Recall that
\begin{align*}
A_{4,0} = h g_{4,1} g_{0,2} , \quad
A_{2,0} = h g_{2,1} g_{0,2} , \quad
A_{2,1} = h g_{2,1} g_{1,2} , \quad
A_{0,2} = h g_{0,1} g_{2,2} , \quad
A_{0,1} = h g_{0,1} g_{2,1} ,
\end{align*}
where $h = 1-2^u = -2^u+1$.

One has
\begin{align*}
x_5 &< 2^u q_2^t  \cdot 10 f_1 t^5 q_1^t
= 2^u 2^{2ut} 2^{2ut^2} \cdot 10 f_1 t^5 \\
&< 2^{u (2t^2 + 4t + 1)} \cdot 10 f_1 t^5
= 2^{u \cdot (2t+1)^2} \cdot 10 f_1 t^5
\end{align*}

Now, $u \cdot (2t+1)^2 < (2t+1)^3$ and $10 x_1 t^5 < 10 t^6 < 10 (2t+1)^6$, so
\begin{align*}
x_5 < 10 \cdot 2^{(2t+1)^3} (2t+1)^6 .
\end{align*}

However, for $z = (2t+1)^3 \geq 3^3 = 27$, one has that $z^2 < 2^z$, so
\begin{align*}
x_5 < 10 \cdot 2^{(2t+1)^3}
= 10 \cdot 2^{(8 x_1 + 3)^3}
= 10 \cdot 2^{512 x_1^3 + 576 x_1^2 + 216 x_1 + 27}.
\end{align*}
\end{proof}

\begin{lemma} \label{proof:f1bounds}
Let $\vec{x} = (x_1,x_2,\ldots,x_8) \in \mathbb{N}^8$. The relation
\begin{align*}
\hat{F}(n,\vec{x},[34]) = 0
& \iff \exists \vec{y} \in \mathbb{N}^{34} :
E_{!}(x_1, [13], x_2)
+ (x_3 - x_2^2)^2
+ (x_4 - x_2 x_3)^2 \\
&\quad + (x_5 - 4 x_2 - 1)^2
+ E_{M4}(x_2, x_3, x_4, x_5, [9], x_6) \\
&\quad + E_{\HW}(x_6, [12], x_7)
+ \left(x_7 + (4 x_2 + 5) \left( -x_5^2 + x_8 - 2^n\right) \right)^2
= 0 ,
\end{align*}
has the bounds:
{\fontsize{8}{12}\selectfont
\begin{align*}
\begin{array}{lll}
y_1 < x_1^2 + 1 ,
& y_2 < 2^{3 y_1 + 1} ,
& y_3 < x_1^3 + 1 , \\
y_4 < 2^{2 y_2 + 4},
& y_5 < 28 y_2^3 + 9,
& y_6 < 10 y_2^2 + 9, \\
y_7 < 2^{28 y_2^3 + 9},
& y_8 < 2^{10 y_2^2 + 8} ,
& y_9 < 2^{10 y_2^2 + 8} , \\
y_{10} < 2^{28 y_2^3 + 9} ,
& y_{11} < 2^{2 y_2 + 4} ,
& y_{12} < 2^{2 y_2 + 4} , \\
y_{13} < 2^{2 y_2 + 4} ,
& y_{14} < 2^{8 x_2 + 10} ,
& y_{15} < 2^{32 x_2 + 48 x_2 + 10} , \\
y_{16} < (4 x_2 + 1) 2^{32 x_2^2 + 48 x_2 + 10} ,
& y_{17} < (4 x_2 + 1)^3 2^{32 x_2^2 + 48 x_2 + 10} ,
& y_{18} < (4 x_2 + 1)^5 2^{32 x_2^2 + 48 x_2 + 10} ,  \\
y_{19} < 2^{128 x_2^3 + 224 x_2^2 + 88 x_2 + 10},
& y_{20} < (4 x_2 + 1)^2 2^{128 x_2^3 + 224 x_2^2 + 88 x_2 + 10},
& y_{21} < (4 x_2 + 1)^3 2^{128 x_2^3 + 224 x_2^2 + 88 x_2 + 10},  \\
y_{22} < 2^{128 x_2^3 + 224 x_2^2 + 92 x_2 + 15} ,
& y_{23} < 2^{4 x_6 + 4} ,
& y_{24} < 224 x_6^3 + 9 ,  \\
y_{25} < 40 x_6^2 + 9 ,
& y_{26} < 2^{224 x_6^3 + 9} ,
& y_{27} < 2^{40 x_6^2 + 8} ,  \\
y_{28} < 2^{40 x_6^2 + 8} ,
& y_{29} < 2^{224 x_6^3 + 9} ,
& y_{30} < 2^{4 x_6 + 4} ,  \\
y_{31} < 2^{4 x_6 + 4} ,
& y_{32} < 2^{4 x_6 + 5} ,
& y_{33} < 2^{4 x_6 + 5} ,  \\
y_{34} < 2^{4 x_6 + 4} + 1 ,
& x_1 < n^2 + 1 ,
& x_2 < 2^{3 y_3 + 1} ,  \\
x_3 < 2^{2(3 x_1^3 + 1)} ,
& x_4 < 2^{3(3 x_1^3 + 1)} ,
& x_5 < 2^{3 x_1^3 + 4} ,  \\
x_6 < 2^{512 x_2^3 + 576 x_2^2 + 216 x_2 + 27} ,
& x_7 < x_6 + 1 ,
& x_8 < 2^{n+2} .
\end{array}
\end{align*}
}
\end{lemma}
\begin{proof}
We start with
\begin{align*}
x_1 < n^2 + 1 .
\end{align*}

By \cref{proof:factorialsingle-foldbounds}, the relation
\begin{align*}
E_{!}(x_1, [13], x_2) = 0 ,
\end{align*}
has the bounds:
\begin{align*}
\begin{array}{lll}
y_1 < x_1^2 + 1 ,
& y_2 < 2^{3 y_1 + 1} ,
& y_3 < x_1^3 + 1 , \\
y_4 < 2^{2 y_2 + 4},
& y_5 < 28 y_2^3 + 9,
& y_6 < 10 y_2^2 + 9, \\
y_7 < 2^{28 y_2^3 + 9},
& y_8 < 2^{10 y_2^2 + 8} ,
& y_9 < 2^{10 y_2^2 + 8} , \\
y_{10} < 2^{28 y_2^3 + 9} ,
& y_{11} < 2^{2 y_2 + 4} ,
& y_{12} < 2^{2 y_2 + 4}  , \\
y_{13} < 2^{2 y_2 + 4}  ,
& x_2 < 2^{3 y_3 + 1} .
\end{array}
\end{align*}

From the above bounds, we obtain
\begin{align*}
x_3 &= x_2^2 < 2^{2(3 x_1^3 + 1)} , \\
x_4 &= x_2^3 < 2^{3(3 x_1^3 + 1)} , \\
x_5 &= 4 x_2 + 1 < 4 \cdot 2^{3 x_1^3 + 1} + 1 < 2^{3 x_1^3 + 4} .
\end{align*}

By \cref{proof:msingle-foldbounds}, the relation
\begin{align*}
E_{M4}(x_2, x_3, x_4, x_5, [9], x_6) = 0 ,
\end{align*}
has the bounds:
\begin{align*}
y_{14} &< 2^{8 x_2 + 10} , \\
y_{15} &< 2^{32 x_2 + 48 x_2 + 10} , \\
y_{16} &< (4 x_2 + 1) 2^{32 x_2^2 + 48 x_2 + 10} , \\
y_{17} &< (4 x_2 + 1)^3 2^{32 x_2^2 + 48 x_2 + 10} , \\
y_{18} &< (4 x_2 + 1)^5 2^{32 x_2^2 + 48 x_2 + 10} , \\
y_{19} &< 2^{128 x_2^3 + 224 x_2^2 + 88 x_2 + 10}, \\
y_{20} &< (4 x_2 + 1)^2 2^{128 x_2^3 + 224 x_2^2 + 88 x_2 + 10}, \\
y_{21} &< (4 x_2 + 1)^3 2^{128 x_2^3 + 224 x_2^2 + 88 x_2 + 10}, \\
y_{22} &< 2^{128 x_2^3 + 224 x_2^2 + 92 x_2 + 15} , \\
x_6 &< 2^{512 x_2^3 + 576 x_2^2 + 216 x_2 + 27} .
\end{align*}

By \cref{proof:hwsingle-foldbounds}, the relation
\begin{align*}
E_{\HW}(x_6, [12], x_7) = 0 ,
\end{align*}
has the bounds:
\begin{align*}
\begin{array}{lll}
y_{23} < 2^{4 x_6 + 4} ,
& y_{24} < 224 x_6^3 + 9 ,
& y_{25} < 40 x_6^2 + 9 , \\
y_{26} < 2^{224 x_6^3 + 9} ,
& y_{27} < 2^{40 x_6^2 + 8} ,
& y_{28} < 2^{40 x_6^2 + 8} , \\
y_{29} < 2^{224 x_6^3 + 9} ,
& y_{30} < 2^{4 x_6 + 4} ,
& y_{31} < 2^{4 x_6 + 4} , \\
y_{32} < 2^{4 x_6 + 5},
& y_{33} < 2^{4 x_6 + 5},
& y_{34} < 2^{4 x_6 + 4} + 1 , \\
x_7 < x_6 + 1 .
\end{array}
\end{align*}

The final relation
\begin{align*}
\left(x_7 + (4 x_2 + 5) \left( -x_5^2 + x_8 - 2^n\right) \right)^2 = 0 ,
\end{align*}
has the bound:
\begin{align*}
x_8 < 2^{n+2} .
\end{align*}
\end{proof}

\begin{lemma} \label{proof:f1upperbound}
Let $\vec{x} = (x_1,x_2,\ldots,x_8)$ and let $\vec{y} = (y_1,y_2,\ldots,y_{34})$. Then
\begin{align*}
\forall (\vec{x}, \vec{y}) \in \mathbb{N}^{42} : \hat{F}(n,\vec{x},[34])=0, \quad \| (\vec{x}, \vec{y}) \|_{\infty} < 2^{2^{2 n^4 + 16}} .
\end{align*}
\end{lemma}
\begin{proof}
From \cref{proof:f1bounds}, it is clear that
\begin{align*}
\| (\vec{x}, \vec{y}) \|_{\infty} \leq 2^{224 \cdot M(4 x_2) + 9} < 2^{224 \cdot 10 \cdot 2^{(8 x_2 + 3)^3}} .
\end{align*}

By \cref{proof:primeinequality}, the largest solution we must consider for a given $n$ is $x_1 = n^2$. Hence, the largest $x_2$ we must consider is $x_2 = (n^2)!$. Making the substitution, we obtain
\begin{align*}
\| (\vec{x}, \vec{y}) \|_{\infty} < 2^{2240 \cdot 2^{(8 \cdot (n^2)! + 3)^3}} .
\end{align*}

Since $(n!)^2 \leq 2^{n^3}$ for all $n \in \mathbb{N}$, we have
\begin{align*}
\| (\vec{x}, \vec{y}) \|_{\infty}
& < 2^{2240 \cdot 2^{(8 \cdot 2^{n^3} + 3)^3}}
= 2^{2240 \cdot 2^{(2^{n^3 + 3} + 3)^3}} \\
& < 2^{2^{12} \cdot 2^{(2^{n^3 + 3} + 3)^3}}
= 2^{2^{(2^{n^3 + 3} + 3)^3 + 12}} .
\end{align*}

Finally, $(2^{n^3 + 3} + 3)^3 + 12 < 2^{2n^4+16}$ for all $n \in \mathbb{N}$, so
\begin{align*}
\| (\vec{x}, \vec{y}) \|_{\infty} < 2^{2^{(2^{n^3 + 3} + 3)^3 + 12}} < 2^{2^{2^{2n^4+16}}} .
\end{align*}
\end{proof}

\section[The n-th prime function]{The $n$-th prime function} \label{section:nthprimeterm}
With our exponential Diophantine equation $\hat{F}(n,\vec{x},[34]) = 0$ developed in \cref{section:monomials} and its bounds established in \cref{section:single-foldbounds}, we are finally ready to write down an explicit arithmetic term for the $n$-th prime number.

Since the hypercube method is typically defined to accept a single tuple of variables $\vec{x}$, we will rewrite $\hat{F}(n,\vec{x},[34]) = 0$ as $\hat{F}(n,\vec{x}) = 0$ by extending the tuple $\vec{x} = (x_1,x_2,\ldots,x_{42})$ and then relabeling the variables $\vec{y} = (y_1,y_2,\ldots,y_{34})$ over the appended $x_i$ variables $(x_{9},x_{10},\ldots,x_{42})$ respectively. Thus,
\begin{align*}
     & \hat{F}(n,\vec{x}) = \hat{F}(n,x_1,x_2,\ldots,x_{42}) = 0 \\
     & \iff
     \hat{F}(n,x_1,x_2,\ldots,x_8,[34]) = \hat{F}(n,x_1,x_2,\ldots,x_8,y_1,y_2,\ldots,y_{34}) = 0 .
\end{align*}
The full monomial expansion of $\hat{F}(n,\vec{x})$ is written in Appendix \cref{appendix:F1monomials}.

\begin{theorem} \label{NthPrimeTerm}
Let $\vec{x} = (x_1,x_2,\ldots,x_{42})$. For all $n \in \Z^+$, the $n$-th prime number, $p(n)$, is given by the arithmetic term:
\begin{align*}
p(n) = \frac{\HW(\hat{Q}(n))}{u(n)} - t(n)^{42} ,
\end{align*}
where $k=42$, $t(n)=2^{2^{2^{2n^4+16}}}$, $u(n)=2^{2^{9 t(n) + 8} + 9}$, and $\hat{Q}(n)$ is the arithmetic term constructed by the hypercube method from the monomial expansion of the exponential Diophantine equation $\hat{F}(n,\vec{x})=0$ (given in Appendix \cref{appendix:F1monomials}) with the parameters $k$, $t(n)$, and $u(n)$.
\end{theorem}
\begin{proof}
By \cref{proof:f1upperbound}, we have that all quantified variables in the solutions to $\hat{F}(n,\vec{x})=0$ are bounded above by
\begin{align*}
t(n) = 2^{2^{2^{2n^4+16}}} .
\end{align*}

The number of monomials in the expansion of $\hat{F}(n,\vec{x})=0$ is
\begin{align*}
\ell = 498 .
\end{align*}

Replacing all variables in $\hat{F}(n,\vec{x}) = 0$ with $t(n)$, we define the exponential polynomial expression
\begin{align*}
\Lambda(n) = \hat{F}(t(n),\ldots,t(n)) .
\end{align*}

Let $\Upsilon(n)$ be the largest monomial in $\Lambda(n)$. To apply the hypercube method, it suffices to set
\begin{align*}
u(n) = \ell \cdot \Upsilon(n) .
\end{align*}

Applying \cref{proof:f1bounds} and \cref{proof:factorialsingle-foldbounds}, we see that
\begin{align*}
\Upsilon(n) &< \ell \cdot 2^{28 \cdot 2^{3(3 t(n) + 1)}}
< \ell \cdot 2^{2^5 \cdot 2^{3(3 t(n) + 1)}}
< \ell \cdot 2^{2^{3(3 t(n) + 1) + 5}} \\
&< \ell \cdot 2^{2^{9 t(n) + 8}}
= 498 \cdot 2^{2^{9 t(n) + 8}}
< 2^9 \cdot 2^{2^{9 t(n) + 8}}
= 2^{2^{9 t(n) + 8} + 9} .
\end{align*}

Hence, we can use
\begin{align*}
u(n) = 2^{2^{9 t(n) + 8} + 9}  .
\end{align*}

Finally, applying the hypercube method as described in \cref{section:hypercube}, it follows from \cref{proof:primeeq2} that
\begin{align*}
\frac{\HW(\hat{Q}(n))}{u(n)} - t(n)^{42}
= | \vec{x} \in \mathbb{N}^{42} : \hat{F}(n,\vec{x})=0 \} |
= p(n) .
\end{align*}
\end{proof}

There is also a somewhat surprising consequence. Consider the term:
\begin{align} \label{TermT}
T(x) = p(\pi(x)+1) ,
\end{align}
where $p(a)$ and $\pi(b)$ are the already built closed terms representing the $a$-th prime and the number of primes which are less or equal to $b$.

\begin{theorem}
The closed term $T(x)$ has the property that for every $x \in \mathbb{N}$, $T(x)$ is the smallest prime strictly bigger than $x$. The recurrent sequence $x(n+1) = T(x(n))$ starting with $x(1)=2$ coincides with the sequence of prime numbers: $\forall n \geq 1, x(n) = p(n)$.
\end{theorem}

So the sequence of prime numbers proves to be a simple recurrent sequence, where the word \textit{simple} means nothing but the fact that any prime depends recurrently only on its predecessor in this sequence.

\begin{appendices}

\clearpage

\section[Monomials for the n-th prime equation \hat{F}=0]{Monomials for the $n$-th prime equation $\hat{F}(n,\vec{x})=0$} \label{appendix:F1monomials}

{\fontsize{8}{12}\selectfont
\begin{align*}
\begin{array}{ll}
\hat{F}(n,\vec{x}) = -2^{184 x_{2} + 288 x_{3} + 128 x_{4} + 35} x_{2}^{2} x_{23}^{3} x_{29}&+ 2^{328 x_{2} + 384 x_{3} + 128 x_{4} + 66} x_{2}^{2} x_{23}^{3} x_{29} \\
-2^{88 x_{2} + 224 x_{3} + 128 x_{4} + 11} x_{23}^{3} x_{29} x_{5}^{2}&-2^{328 x_{2} + 384 x_{3} + 128 x_{4} + 63} x_{2} x_{23}^{2} x_{28} \\
+ 2^{328 x_{2} + 384 x_{3} + 128 x_{4} + 64} x_{2} x_{23}^{3} x_{29}&-3 \cdot 2^{32 x_{2} + 32 x_{3} + 40} x_{2}^{4} x_{22}^{5} x_{26} \\
-3 \cdot 2^{40 x_{2} + 32 x_{3} + 46} x_{2}^{2} x_{22}^{5} x_{26}&-3 \cdot 2^{32 x_{2} + 32 x_{3} + 39} x_{2}^{3} x_{22}^{5} x_{26} \\
+ 2^{88 x_{2} + 224 x_{3} + 128 x_{4} + 11} x_{23}^{2} x_{28} x_{5}&+ 3 \cdot 2^{56 x_{2} + 32 x_{3} + 29} x_{2}^{3} x_{22}^{5} x_{26} \\
+ 3 \cdot 2^{56 x_{2} + 32 x_{3} + 26} x_{2}^{2} x_{22}^{5} x_{26}&+ 3 \cdot 2^{32 x_{2} + 32 x_{3} + 36} x_{2}^{2} x_{22}^{5} x_{26} \\
-2^{416 x_{2} + 608 x_{3} + 256 x_{4} + 76} x_{2}^{2} x_{5}^{2}&+ 2^{272 x_{2} + 512 x_{3} + 256 x_{4} + 45} x_{2}^{2} x_{5}^{2} \\
-3 \cdot 2^{56 x_{2} + 32 x_{3} + 25} x_{2} x_{22}^{5} x_{26}&-2^{328 x_{2} + 384 x_{3} + 128 x_{4} + 61} x_{23}^{3} x_{29} \\
+ 3 \cdot 2^{32 x_{2} + 32 x_{3} + 35} x_{2} x_{22}^{5} x_{26}&-2^{416 x_{2} + 608 x_{3} + 256 x_{4} + 74} x_{2} x_{5}^{2} \\
-1649267441647 \cdot 2^{72 x_{2} + 32 x_{3} + 47} x_{2}^{2}&+ 1099511627787 \cdot 2^{72 x_{2} + 32 x_{3} + 51} x_{2}^{4} \\
+ 1099511627809 \cdot 2^{72 x_{2} + 32 x_{3} + 49} x_{2}^{3}&-2^{16 x_{2} + 32 x_{3} + 19} x_{2}^{4} x_{22}^{5} x_{26} \\
-2^{16 x_{2} + 32 x_{3} + 11} x_{22}^{5} x_{26} x_{5}^{4}&-3 \cdot 2^{72 x_{2} + 64 x_{3} + 39} x_{2}^{3} x_{5}^{4} \\
-3 \cdot 2^{72 x_{2} + 64 x_{3} + 36} x_{2}^{2} x_{5}^{4}&-3 \cdot 2^{48 x_{2} + 64 x_{3} + 46} x_{2}^{2} x_{5}^{4} \\
-2^{32 x_{2} + 32 x_{3} + 35} x_{2}^{2} x_{22}^{3} x_{25}&-2^{16 x_{2} + 32 x_{3} + 11} x_{22}^{3} x_{25} x_{5}^{2} \\
-2^{280 x_{2} + 352 x_{3} + 128 x_{4} + 51} x_{23} x_{27}&+ 2^{40 x_{2} + 32 x_{3} + 51} x_{2}^{4} x_{22}^{5} x_{26} \\
+ 2^{56 x_{2} + 32 x_{3} + 31} x_{2}^{4} x_{22}^{5} x_{26}&+ 2^{40 x_{2} + 32 x_{3} + 49} x_{2}^{3} x_{22}^{5} x_{26} \\
+ 3 \cdot 2^{48 x_{2} + 64 x_{3} + 50} x_{2}^{4} x_{5}^{4}&+ 3 \cdot 2^{48 x_{2} + 64 x_{3} + 49} x_{2}^{3} x_{5}^{4} \\
+ 2^{56 x_{2} + 32 x_{3} + 26} x_{2}^{2} x_{22}^{3} x_{25}&-11 \cdot 2^{56 x_{2} + 32 x_{3} + 21} x_{22}^{5} x_{26} \\
+ 3 \cdot 2^{56 x_{2} + 64 x_{3} + 56} x_{2}^{2} x_{5}^{4}&-11 \cdot 2^{32 x_{2} + 32 x_{3} + 31} x_{22}^{5} x_{26} \\
-2^{416 x_{2} + 608 x_{3} + 256 x_{4} + 73} x_{2} x_{5}&-2^{8 x_{2} + 18} x_{2}^{3} x_{24} x_{25} x_{28} x_{29} \\
-2^{8 x_{2} + 16} x_{2}^{2} x_{24} x_{25} x_{27} x_{29}&-2^{4 x_{2} + 11} x_{2}^{2} x_{24} x_{26} x_{27} x_{29} \\
-2^{512 x_{2} + 672 x_{3} + 256 x_{4} + 100} x_{2}^{4}&+ 2^{4 x_{2} + 14} x_{2}^{3} x_{24} x_{25} x_{28} x_{29} \\
+ 2^{4 x_{2} + 12} x_{2}^{2} x_{24} x_{25} x_{27} x_{29}&+ 2^{8 x_{2} + 15} x_{2}^{2} x_{24} x_{26} x_{27} x_{29} \\
+ 2199023255487 \cdot 2^{72 x_{2} + 32 x_{3} + 44} x_{2}&-2^{512 x_{2} + 672 x_{3} + 256 x_{4} + 98} x_{2}^{3} \\
-2^{280 x_{2} + 352 x_{3} + 128 x_{4} + 55} x_{2}^{2}&+ 2^{656 x_{2} + 768 x_{3} + 256 x_{4} + 130} x_{2}^{4} \\
-2^{184 x_{2} + 288 x_{3} + 128 x_{4} + 31} x_{5}^{2}&+ 2^{656 x_{2} + 768 x_{3} + 256 x_{4} + 129} x_{2}^{3} \\
-3 \cdot 2^{48 x_{2} + 64 x_{3} + 45} x_{2} x_{5}^{4}&+ 2^{656 x_{2} + 768 x_{3} + 256 x_{4} + 124} x_{2}^{2} \\
-2^{232 x_{2} + 320 x_{3} + 128 x_{4} + 45} x_{2}^{2}&-2^{136 x_{2} + 256 x_{3} + 128 x_{4} + 21} x_{5}^{2} \\
+ 2^{368 x_{2} + 576 x_{3} + 256 x_{4} + 68} x_{2}^{4}&+ 2^{40 x_{2} + 32 x_{3} + 45} x_{2} x_{22}^{5} x_{26} \\
+ 2^{512 x_{2} + 672 x_{3} + 256 x_{4} + 95} x_{2}^{2}&+ 2^{424 x_{2} + 448 x_{3} + 128 x_{4} + 86} x_{2}^{2} \\
+ 2^{416 x_{2} + 608 x_{3} + 256 x_{4} + 71} x_{5}^{2}&+ 2^{176 x_{2} + 448 x_{3} + 256 x_{4} + 20} x_{5}^{4} \\
+ 3 \cdot 2^{72 x_{2} + 64 x_{3} + 35} x_{2} x_{5}^{4}&+ 2^{56 x_{2} + 32 x_{3} + 24} x_{2} x_{22}^{3} x_{25} \\
-2^{4 x_{2} + 14} x_{2}^{3} x_{24}^{2} x_{25} x_{29}&+ 2^{376 x_{2} + 416 x_{3} + 128 x_{4} + 76} x_{2}^{2} \\
-2^{8 x_{2} + 17} x_{2}^{2} x_{24} x_{25}^{2} x_{28}&+ 2^{176 x_{2} + 448 x_{3} + 256 x_{4} + 20} x_{5}^{2} \\
+ 2^{8 x_{2} + 18} x_{2}^{3} x_{24}^{2} x_{25} x_{29}&+ 2^{4 x_{2} + 13} x_{2}^{2} x_{24} x_{25}^{2} x_{28} \\
-2^{4 x_{2} + 10} x_{2} x_{24} x_{25} x_{26} x_{27}&-2^{8 x_{2} + 14} x_{2} x_{25} x_{26} x_{27} x_{28} \\
-2^{656 x_{2} + 768 x_{3} + 256 x_{4} + 124} x_{2}&+ 2^{8 x_{2} + 14} x_{2} x_{24} x_{25} x_{26} x_{27} \\
+ 2^{4 x_{2} + 10} x_{2} x_{25} x_{26} x_{27} x_{28}&-2^{288 x_{2} + 352 x_{3} + 128 x_{4} + 61} x_{30} \\
-2^{284 x_{2} + 352 x_{3} + 128 x_{4} + 56} x_{30}&-2^{72 x_{2} + 64 x_{3} + 41} x_{2}^{4} x_{5}^{4} \\
-2^{56 x_{2} + 64 x_{3} + 61} x_{2}^{4} x_{5}^{4}&-2^{56 x_{2} + 64 x_{3} + 59} x_{2}^{3} x_{5}^{4} \\
-2^{4 x_{2} + 14} x_{2}^{4} x_{24}^{2} x_{29}^{2}&-2^{72 x_{2} + 64 x_{3} + 36} x_{2}^{2} x_{5}^{2} \\
-21 \cdot 2^{112 x_{2} + 64 x_{3} + 52} x_{2}^{4}&-2^{4 x_{2} + 12} x_{2}^{2} x_{24}^{2} x_{25}^{2} \\
-2^{4 x_{2} + 12} x_{2}^{2} x_{25}^{2} x_{28}^{2}&-275 \cdot 2^{88 x_{2} + 64 x_{3} + 55} x_{2}^{2} \\
-25 \cdot 2^{112 x_{2} + 64 x_{3} + 50} x_{2}^{3}&-119 \cdot 2^{64 x_{2} + 64 x_{3} + 68} x_{2}^{4} \\
-15 \cdot 2^{88 x_{2} + 64 x_{3} + 68} x_{2}^{7}&-39 \cdot 2^{88 x_{2} + 64 x_{3} + 65} x_{2}^{6} \\
+ 2^{32 x_{2} + 64 x_{3} + 29} x_{2}^{4} x_{5}^{4}&-21 \cdot 2^{96 x_{2} + 64 x_{3} + 71} x_{2}^{4} \\
+ 2^{8 x_{2} + 18} x_{2}^{4} x_{24}^{2} x_{29}^{2}&+ 2^{424 x_{2} + 448 x_{3} + 128 x_{4} + 84} x_{2} \\
-37 \cdot 2^{72 x_{2} + 64 x_{3} + 80} x_{2}^{4}&-33 \cdot 2^{56 x_{2} + 32 x_{3} + 70} x_{2}^{4} \\
-13 \cdot 2^{48 x_{2} + 32 x_{3} + 59} x_{2}^{4}&-33 \cdot 2^{56 x_{2} + 32 x_{3} + 69} x_{2}^{3} \\
-33 \cdot 2^{48 x_{2} + 32 x_{3} + 60} x_{2}^{4}&+ 2^{8 x_{2} + 16} x_{2}^{2} x_{24}^{2} x_{25}^{2} \\
+ 2^{8 x_{2} + 16} x_{2}^{2} x_{25}^{2} x_{28}^{2}&+ 2^{48 x_{2} + 64 x_{3} + 45} x_{2}^{2} x_{5}^{2} \\
-11 \cdot 2^{40 x_{2} + 32 x_{3} + 41} x_{5}^{4}&+ 2^{376 x_{2} + 416 x_{3} + 128 x_{4} + 74} x_{2} \\
-15 \cdot 2^{64 x_{2} + 32 x_{3} + 77} x_{2}^{2}&-33 \cdot 2^{48 x_{2} + 32 x_{3} + 59} x_{2}^{3} \\
-11 \cdot 2^{32 x_{2} + 32 x_{3} + 39} x_{2}^{4}&-2^{8 x_{2} + 15} x_{2} x_{24} x_{25}^{2} x_{27}+
\end{array}
\end{align*}
\begin{align*}
\begin{array}{lll}
-2^{4 x_{2} + 11} x_{2} x_{25}^{2} x_{27} x_{28}&-2^{4 x_{2} + 10} x_{2}^{2} x_{24} x_{29} x_{30}&-11 \cdot 2^{32 x_{2} + 32 x_{3} + 31} x_{5}^{4} \\
-33 \cdot 2^{56 x_{2} + 32 x_{3} + 66} x_{2}^{2}&-19 \cdot 2^{40 x_{2} + 32 x_{3} + 45} x_{2}^{2}&-549755813949 \cdot 2^{72 x_{2} + 32 x_{3} + 42} \\
-3 \cdot 2^{72 x_{2} + 64 x_{3} + 90} x_{2}^{8}&-3 \cdot 2^{88 x_{2} + 64 x_{3} + 70} x_{2}^{8}&-9 \cdot 2^{72 x_{2} + 64 x_{3} + 88} x_{2}^{7} \\
-3 \cdot 2^{72 x_{2} + 64 x_{3} + 47} x_{2}^{7}&-5 \cdot 2^{96 x_{2} + 64 x_{3} + 75} x_{2}^{5}&-3 \cdot 2^{72 x_{2} + 64 x_{3} + 44} x_{2}^{6} \\
+ 3 \cdot 2^{112 x_{2} + 64 x_{3} + 59} x_{2}^{7}&+ 15 \cdot 2^{72 x_{2} + 64 x_{3} + 84} x_{2}^{5}&+ 15 \cdot 2^{88 x_{2} + 64 x_{3} + 64} x_{2}^{5} \\
+ 3 \cdot 2^{112 x_{2} + 64 x_{3} + 58} x_{2}^{6}&-3 \cdot 2^{48 x_{2} + 64 x_{3} + 54} x_{2}^{6}&-2^{40 x_{2} + 32 x_{3} + 41} x_{22}^{5} x_{26} \\
+ 43 \cdot 2^{88 x_{2} + 64 x_{3} + 62} x_{2}^{4}&+ 11 \cdot 2^{72 x_{2} + 64 x_{3} + 39} x_{2}^{4}&+ 19 \cdot 2^{64 x_{2} + 32 x_{3} + 80} x_{2}^{4} \\
+ 3 \cdot 2^{112 x_{2} + 64 x_{3} + 54} x_{2}^{5}&-3 \cdot 2^{48 x_{2} + 64 x_{3} + 53} x_{2}^{5}&+ 11 \cdot 2^{72 x_{2} + 64 x_{3} + 31} x_{5}^{4} \\
-5 \cdot 2^{72 x_{2} + 64 x_{3} + 81} x_{2}^{3}&-3 \cdot 2^{48 x_{2} + 32 x_{3} + 59} x_{2}^{3}&+ 15 \cdot 2^{80 x_{2} + 64 x_{3} + 90} x_{2}^{4} \\
+ 11 \cdot 2^{48 x_{2} + 64 x_{3} + 49} x_{2}^{4}&+ 11 \cdot 2^{80 x_{2} + 32 x_{3} + 61} x_{2}^{4}&+ 11 \cdot 2^{56 x_{2} + 32 x_{3} + 71} x_{2}^{4} \\
-3 \cdot 2^{64 x_{2} + 32 x_{3} + 41} x_{2}^{4}&-9 \cdot 2^{64 x_{2} + 64 x_{3} + 75} x_{2}^{5}&-2^{56 x_{2} + 32 x_{3} + 21} x_{22}^{3} x_{25} \\
+ 11 \cdot 2^{48 x_{2} + 64 x_{3} + 41} x_{5}^{4}&+ 27 \cdot 2^{72 x_{2} + 64 x_{3} + 77} x_{2}^{2}&+ 33 \cdot 2^{56 x_{2} + 32 x_{3} + 66} x_{2}^{2} \\
-7 \cdot 2^{80 x_{2} + 64 x_{3} + 89} x_{2}^{3}&+ 33 \cdot 2^{80 x_{2} + 32 x_{3} + 59} x_{2}^{3}&+ 11 \cdot 2^{56 x_{2} + 32 x_{3} + 69} x_{2}^{3} \\
-9 \cdot 2^{64 x_{2} + 32 x_{3} + 39} x_{2}^{3}&+ 2^{4 x_{2} + 11} x_{2} x_{24} x_{25}^{2} x_{27}&+ 2^{8 x_{2} + 15} x_{2} x_{25}^{2} x_{27} x_{28} \\
+ 3 \cdot 2^{112 x_{2} + 64 x_{3} + 46} x_{2}^{2}&+ 33 \cdot 2^{80 x_{2} + 32 x_{3} + 56} x_{2}^{2}&+ 65 \cdot 2^{48 x_{2} + 32 x_{3} + 55} x_{2}^{2} \\
+ 51 \cdot 2^{64 x_{2} + 64 x_{3} + 69} x_{2}^{3}&+ 3 \cdot 2^{96 x_{2} + 64 x_{3} + 77} x_{2}^{6}&+ 3 \cdot 2^{48 x_{2} + 64 x_{3} + 58} x_{2}^{8} \\
+ 3 \cdot 2^{72 x_{2} + 64 x_{3} + 85} x_{2}^{6}&+ 3 \cdot 2^{48 x_{2} + 64 x_{3} + 57} x_{2}^{7}&+ 9 \cdot 2^{64 x_{2} + 64 x_{3} + 78} x_{2}^{8} \\
+ 3 \cdot 2^{72 x_{2} + 64 x_{3} + 43} x_{2}^{5}&+ 3 \cdot 2^{56 x_{2} + 64 x_{3} + 64} x_{2}^{6}&+ 9 \cdot 2^{64 x_{2} + 64 x_{3} + 78} x_{2}^{7} \\
+ 5 \cdot 2^{96 x_{2} + 64 x_{3} + 70} x_{2}^{3}&+ 9 \cdot 2^{64 x_{2} + 64 x_{3} + 75} x_{2}^{6}&+ 3 \cdot 2^{96 x_{2} + 64 x_{3} + 67} x_{2}^{2} \\
+ 3 \cdot 2^{88 x_{2} + 32 x_{3} + 69} x_{2}^{3}&+ 7 \cdot 2^{40 x_{2} + 32 x_{3} + 49} x_{2}^{4}&+ 3 \cdot 2^{24 x_{2} + 32 x_{3} + 29} x_{2}^{4} \\
+ 3 \cdot 2^{24 x_{2} + 32 x_{3} + 21} x_{5}^{4}&+ 3 \cdot 2^{88 x_{2} + 32 x_{3} + 66} x_{2}^{2}&+ 9 \cdot 2^{48 x_{2} + 32 x_{3} + 56} x_{2}^{2} \\
+ 9 \cdot 2^{40 x_{2} + 32 x_{3} + 49} x_{2}^{3}&+ 2^{96 x_{2} + 64 x_{3} + 21} x_{23}^{3} x_{29}&+ 2^{48 x_{2} + 32 x_{3} + 11} x_{23}^{3} x_{29} \\
+ 2^{4 x_{2} + 10} x_{2}^{2} x_{24} x_{29} x_{6}&+ 7 \cdot 2^{80 x_{2} + 64 x_{3} + 86} x_{2}^{2}&+ 3 \cdot 2^{64 x_{2} + 64 x_{3} + 66} x_{2}^{2} \\
-2^{56 x_{2} + 64 x_{3} + 55} x_{2} x_{5}^{4}&-2^{72 x_{2} + 64 x_{3} + 34} x_{2} x_{5}^{2}&+ 11 \cdot 2^{24 x_{2} + 31} x_{22}^{5} x_{26} \\
+ 11 \cdot 2^{16 x_{2} + 21} x_{22}^{5} x_{26}&+ 65 \cdot 2^{112 x_{2} + 64 x_{3} + 44} x_{2}&-33 \cdot 2^{80 x_{2} + 32 x_{3} + 55} x_{2} \\
-33 \cdot 2^{64 x_{2} + 64 x_{3} + 65} x_{2}&-2^{424 x_{2} + 448 x_{3} + 128 x_{4} + 81}&-3 \cdot 2^{8 x_{2} + 11} x_{22}^{5} x_{26} \\
+ 2^{656 x_{2} + 768 x_{3} + 256 x_{4} + 120}&-2^{376 x_{2} + 416 x_{3} + 128 x_{4} + 71}&-7 \cdot 2^{72 x_{2} + 64 x_{3} + 76} x_{2} \\
-3 \cdot 2^{88 x_{2} + 32 x_{3} + 65} x_{2}&+ 33 \cdot 2^{56 x_{2} + 32 x_{3} + 65} x_{2}&-3 \cdot 2^{48 x_{2} + 32 x_{3} + 55} x_{2} \\
-2^{4 x_{2} + 9} x_{2} x_{24} x_{25} x_{30}&+ 2^{568 x_{2} + 704 x_{3} + 256 x_{4} + 110}&-2^{288 x_{2} + 352 x_{3} + 128 x_{4} + 61} \\
+ 11 \cdot 2^{56 x_{2} + 32 x_{3} + 65} x_{2}&+ 19 \cdot 2^{64 x_{2} + 32 x_{3} + 34} x_{2}&+ 33 \cdot 2^{48 x_{2} + 32 x_{3} + 55} x_{2} \\
-9 \cdot 2^{40 x_{2} + 32 x_{3} + 45} x_{2}&+ 2^{560 x_{2} + 704 x_{3} + 256 x_{4} + 100}&-2^{280 x_{2} + 352 x_{3} + 128 x_{4} + 51} \\
-256 x_{2}^{3} x_{24} x_{25} x_{28} x_{29}&+ 7 \cdot 2^{64 x_{2} + 32 x_{3} + 76} x_{2}&-2^{8 x_{2} + 12} x_{25} x_{26} x_{27}^{2} \\
+ 2^{4 x_{2} + 9} x_{2} x_{25} x_{28} x_{30}&-2^{4 x_{2} + 9} x_{2} x_{25} x_{28} x_{6}&-64 x_{2}^{2} x_{24} x_{25} x_{27} x_{29} \\
+ 2^{4 x_{2} + 9} x_{2} x_{24} x_{25} x_{6}&+ 32 x_{2}^{2} x_{24} x_{26} x_{27} x_{29}&+ 2^{4 x_{2} + 8} x_{25} x_{26} x_{27}^{2} \\
-2^{72 x_{2} + 64 x_{3} + 49} x_{2}^{8}&+ 2^{80 x_{2} + 64 x_{3} + 100} x_{2}^{8}&-2^{56 x_{2} + 64 x_{3} + 69} x_{2}^{8} \\
+ 2^{112 x_{2} + 64 x_{3} + 60} x_{2}^{8}&-2^{56 x_{2} + 64 x_{3} + 67} x_{2}^{7}&-2^{80 x_{2} + 64 x_{3} + 97} x_{2}^{6} \\
-2^{80 x_{2} + 64 x_{3} + 94} x_{2}^{5}&-2^{56 x_{2} + 64 x_{3} + 63} x_{2}^{5}&-2^{48 x_{2} + 32 x_{3} + 51} x_{5}^{4} \\
-2^{88 x_{2} + 64 x_{3} + 58} x_{2}^{3}&-128 x_{2}^{2} x_{24} x_{25}^{2} x_{28}&-2^{64 x_{2} + 32 x_{3} + 39} x_{2}^{2} \\
-2^{24 x_{2} + 32 x_{3} + 21} x_{5}^{2}&-121 \cdot 2^{56 x_{2} + 32 x_{3} + 61}&-2^{32 x_{2} + 32 x_{3} + 31} x_{5}^{2} \\
-121 \cdot 2^{80 x_{2} + 32 x_{3} + 51}&-121 \cdot 2^{48 x_{2} + 32 x_{3} + 51}&+ 2^{96 x_{2} + 64 x_{3} + 81} x_{2}^{8} \\
+ 2^{96 x_{2} + 64 x_{3} + 81} x_{2}^{7}&+ 2^{80 x_{2} + 64 x_{3} + 99} x_{2}^{7}&+ 2^{32 x_{2} + 64 x_{3} + 36} x_{2}^{8} \\
+ 2^{32 x_{2} + 64 x_{3} + 20} x_{5}^{8}&+ 2^{88 x_{2} + 32 x_{3} + 71} x_{2}^{4}&+ 2^{64 x_{2} + 32 x_{3} + 82} x_{2}^{3} \\
+ 2^{56 x_{2} + 64 x_{3} + 59} x_{2}^{4}&+ 256 x_{2}^{3} x_{24}^{2} x_{25} x_{29}&+ 2^{56 x_{2} + 64 x_{3} + 51} x_{5}^{4} \\
+ 2^{72 x_{2} + 64 x_{3} + 31} x_{5}^{2}&+ 2^{32 x_{2} + 64 x_{3} + 20} x_{5}^{4}&+ 2^{8 x_{2} + 12} x_{25}^{2} x_{27}^{2} \\
-2^{4 x_{2} + 8} x_{25}^{2} x_{27}^{2}&+ 2^{8 x_{2} + 10} x_{26}^{2} x_{27}^{2}&-2^{4 x_{2} + 6} x_{26}^{2} x_{27}^{2} \\
+ 121 \cdot 2^{88 x_{2} + 64 x_{3} + 51}&-11 \cdot 2^{88 x_{2} + 32 x_{3} + 61}&-11 \cdot 2^{64 x_{2} + 32 x_{3} + 72} \\
+ 61 \cdot 2^{112 x_{2} + 64 x_{3} + 41}&-11 \cdot 2^{56 x_{2} + 32 x_{3} + 61}&+ 121 \cdot 2^{64 x_{2} + 64 x_{3} + 60} \\
+ 11 \cdot 2^{96 x_{2} + 64 x_{3} + 61}&-16 x_{2} x_{25} x_{26} x_{27} x_{28}&+ 11 \cdot 2^{72 x_{2} + 64 x_{3} + 71} \\
-2^{2 x_{10} + x_{14} + 5} x_{15}^{2}&-2^{4 x_{2} + 6} x_{26} x_{27} x_{30}&+ 33 \cdot 2^{40 x_{2} + 32 x_{3} + 41} \\
+ 16 x_{2} x_{24} x_{25} x_{26} x_{27}&+ 3 \cdot 2^{48 x_{2} + 32 x_{3} + 51}&+ 2^{4 x_{2} + 7} x_{25} x_{27} x_{30} \\
-2^{x_{33} + 4 x_{6} + 5} x_{34}^{2}&-2^{4 x_{2} + 7} x_{25} x_{27} x_{6}&-2^{48 x_{2} + 32 x_{3} + 11} x_{23} \\
+ 2^{32 x_{2} + 41} x_{22}^{5} x_{26}&+ 256 x_{2}^{4} x_{24}^{2} x_{29}^{2}&-2^{96 x_{2} + 64 x_{3} + 68} x_{2} \\
+ 2^{16 x_{2} + 21} x_{22}^{3} x_{25}&-2^{80 x_{2} + 64 x_{3} + 85} x_{2}&+ 2^{4 x_{2} + 6} x_{26} x_{27} x_{6} \\
+ 64 x_{2}^{2} x_{24}^{2} x_{25}^{2}&+ 64 x_{2}^{2} x_{25}^{2} x_{28}^{2}&+ 2^{8 x_{2} + 11} x_{22}^{3} x_{25} \\
-32 x_{2} x_{24} x_{25}^{2} x_{27}&+ 32 x_{2} x_{25}^{2} x_{27} x_{28}&+ 32 x_{2}^{2} x_{24} x_{29} x_{30}+
\end{array}
\end{align*}
\begin{align*}
\begin{array}{lllll}
-32 x_{2}^{2} x_{24} x_{29} x_{6}&+ 5 \cdot 2^{n + 2} x_{2} x_{5}^{2}&+ 2^{2 x_{10} + x_{13} + 5} x_{15}&-2^{x_{13} + x_{14} + 1} x_{15}&-2^{2 x_{10} + 5} x_{15} x_{16} \\
-2^{2 x_{10} + 5} x_{15} x_{18}&-2^{x_{32} + x_{33} + 1} x_{34}&+ 2^{x_{32} + 4 x_{6} + 5} x_{34}&-16 x_{2} x_{25} x_{28} x_{30}&+ 2^{2 x_{10} + 5} x_{12} x_{18} \\
-2^{3 x_{11} + 1} x_{12} x_{2}&-2^{4 x_{6} + 5} x_{34} x_{35}&-2^{4 x_{6} + 5} x_{34} x_{37}&-5 \cdot 2^{n + 2} x_{2} x_{8}&+ 2^{192 x_{2} + 128 x_{3} + 40} \\
+ 2^{n + 1} x_{2}^{2} x_{5}^{2}&+ 2^{144 x_{2} + 96 x_{3} + 31}&+ 16 x_{2} x_{24} x_{25} x_{30}&-16 x_{2} x_{24} x_{25} x_{6}&+ 2^{4 x_{6} + 5} x_{31} x_{37} \\
-2^{48 x_{2} + 32 x_{3} + 11}&-2 x_{2}^{2} x_{5}^{2} x_{8}&+ 16 x_{2} x_{25} x_{28} x_{6}&+ 2^{80 x_{2} + 64 x_{3} + 80}&+ 2^{64 x_{2} + 32 x_{3} + 36} \\
+ 2^{x_{14} + 1} x_{15} x_{16}&+ 2^{x_{33} + 1} x_{34} x_{35}&-2^{x_{7} + 1} x_{31} x_{39}&-2^{x_{7} + 2} x_{31} x_{40}&+ 25 \cdot 2^{n + 1} x_{5}^{2} \\
+ 2^{96 x_{2} + 64 x_{3} + 21}&-4 x_{25} x_{26} x_{27}^{2}&+ 2^{4 x_{10} + 8} x_{15}^{2}&+ 2^{2 x_{10} + 5} x_{15}^{2}&+ 2^{4 x_{10} + 8} x_{18}^{2} \\
+ 2^{8 x_{2} + 10} x_{30}^{2}&+ 2^{x_{7} + 2} x_{40} x_{41}&+ 121 \cdot 2^{48 x_{2} + 60}&+ 121 \cdot 2^{40 x_{2} + 51}&-2^{n + 1} x_{2}^{2} x_{8} \\
+ 11 \cdot 2^{48 x_{2} + 61}&-2^{x_{14} + 1} x_{15}^{2}&+ 2^{4 x_{2} + 6} x_{30}^{2}&-2^{x_{33} + 1} x_{34}^{2}&+ 2^{8 x_{6} + 8} x_{34}^{2} \\
+ 2^{4 x_{6} + 5} x_{34}^{2}&+ 2^{8 x_{6} + 8} x_{37}^{2}&+ 2^{2 x_{7} + 2} x_{40}^{2}&+ 11 \cdot 2^{56 x_{2} + 71}&-3 \cdot 2^{40 x_{2} + 51} \\
+ 11 \cdot 2^{16 x_{2} + 20}&-20 x_{2} x_{5}^{2} x_{8}&+ 5 \cdot 2^{2 n + 1} x_{2}&-25 \cdot 2^{n + 1} x_{8}&+ 7 \cdot 2^{32 x_{2} + 43} \\
-2 x_{2} x_{5}^{2} x_{7}&-2^{2 x_{10} + 5} x_{12}&-2^{2 x_{10} + 5} x_{19}&-2^{3 x_{11} + 1} x_{20}&-2^{8 x_{2} + 11} x_{22} \\
-5 \cdot 2^{n + 1} x_{7}&+ 4 x_{25}^{2} x_{27}^{2}&+ 2^{2 x_{14}} x_{15}^{2}&-2 x_{22} x_{23} x_{24}&-4 x_{25} x_{27} x_{30} \\
+ 2^{2 x_{33}} x_{34}^{2}&-2^{3 x_{1} + 1} x_{10}&+ 2^{8 x_{2} + 11} x_{30}&-2^{4 x_{6} + 5} x_{31}&-2^{4 x_{6} + 5} x_{38} \\
+ x_{22}^{10} x_{26}^{2}&+ 4 x_{1}^{2} x_{10}^{2}&-4 x_{1} x_{10} x_{13}&+ 2 x_{26} x_{27} x_{30}&+ 2^{2 x_{7}} x_{39}^{2} \\
-2 x_{26} x_{27} x_{6}&-2^{n + 1} x_{2} x_{7}&-2^{x_{13} + 1} x_{16}&+ 2^{4 x_{2} + 6} x_{30}&-2^{x_{32} + 1} x_{35} \\
+ x_{22}^{6} x_{25}^{2}&+ x_{23}^{6} x_{29}^{2}&+ x_{23}^{4} x_{28}^{2}&+ x_{22}^{2} x_{24}^{2}&+ x_{23}^{2} x_{27}^{2} \\
+ x_{26}^{2} x_{27}^{2}&-16 x_{10}^{2} x_{13}&+ 2 x_{12} x_{2} x_{20}&+ 4 x_{25} x_{27} x_{6}&-2 x_{1} x_{11} x_{9} \\
+ 2^{x_{13} + 1} x_{15}&+ 2^{x_{32} + 1} x_{34}&-2^{x_{7} + 2} x_{41}&-2^{x_{7} + 2} x_{42}&-4 x_{10}^{3} x_{13} \\
+ x_{12}^{2} x_{2}^{2}&-32 x_{32} x_{6}^{3}&+ 112 x_{1} x_{10}^{2}&-4 x_{10}^{2} x_{14}&-2 x_{2} x_{3} x_{4} \\
-72 x_{32} x_{6}^{2}&-16 x_{33} x_{6}^{2}&+ 2^{x_{7} + 2} x_{40}&+ x_{2}^{2} x_{5}^{4}&+ 48 x_{1} x_{10}^{3} \\
+ x_{2}^{2} x_{3}^{2}&+ x_{2}^{2} x_{8}^{2}&+ x_{1}^{2} x_{9}^{2}&+ 16 x_{1}^{2} x_{10}&-10 x_{5}^{2} x_{7} \\
-50 x_{5}^{2} x_{8}&+ 2 x_{2} x_{7} x_{8}&+ 8 x_{1} x_{10}^{4}&+ 10 x_{2} x_{5}^{4}&-2 x_{2}^{2} x_{3} \\
+ 10 x_{2} x_{8}^{2}&-2 x_{1}^{2} x_{9}&-2^{24 x_{2} + 36}&+ 2^{2 n} x_{2}^{2}&+ 2^{64 x_{2} + 80} \\
-24 x_{10} x_{13}&-16 x_{10} x_{14}&-2^{2 x_{10} + 5}&-2^{8 x_{2} + 10}&+ 128 x_{1} x_{10} \\
-2 x_{12} x_{15}&-4 x_{15} x_{16}&-2 x_{15} x_{17}&-2 x_{31} x_{34}&-4 x_{34} x_{35} \\
-2 x_{34} x_{36}&-2 x_{31} x_{41}&-56 x_{32} x_{6}&-32 x_{33} x_{6}&+ 25 \cdot 2^{2 n} \\
+ 2^{4 x_{10} + 8}&-2^{4 x_{6} + 5}&-8 x_{1} x_{13}&+ 2 x_{16} x_{17}&+ 2 x_{12} x_{19} \\
-2 x_{2} x_{20}&-2 x_{2} x_{21}&+ 2 x_{20} x_{21}&+ 2 x_{23} x_{27}&+ 2 x_{35} x_{36} \\
+ 2 x_{31} x_{38}&+ 2 x_{41} x_{42}&-2 x_{30} x_{6}&+ 2^{8 x_{6} + 8}&+ 2^{2 x_{7} + 2} \\
+ 1152 x_{6}^{5}&+ 116 x_{10}^{4}&+ 2256 x_{6}^{4}&+ 256 x_{10}^{3}&+ 2528 x_{6}^{3} \\
+ 369 x_{10}^{2}&-8 x_{2} x_{5}&+ 1745 x_{6}^{2}&+ 10 x_{7} x_{8}&-2^{x_{7} + 3} \\
+ 256 x_{6}^{6}&+ 32 x_{10}^{5}&+ 4 x_{10}^{6}&+ 25 x_{5}^{4}&+ 16 x_{1}^{2} \\
+ 2 x_{12}^{2}&+ 3 x_{15}^{2}&+ 2 x_{16}^{2}&+ 17 x_{2}^{2}&+ 2 x_{20}^{2} \\
+ 2 x_{23}^{2}&+ 2 x_{30}^{2}&+ 4 x_{31}^{2}&+ 3 x_{34}^{2}&+ 2 x_{35}^{2} \\
+ 2 x_{41}^{2}&+ 25 x_{8}^{2}&+ 2^{6 x_{11}}&+ 2^{2 x_{13}}&+ 2^{2 x_{32}} \\
+ 2^{6 x_{1}}&+ x_{11}^{2}&+ x_{13}^{2}&+ x_{14}^{2}&+ x_{17}^{2} \\
+ x_{19}^{2}&+ x_{21}^{2}&+ x_{22}^{2}&+ x_{32}^{2}&+ x_{33}^{2} \\
+ x_{36}^{2}&+ x_{38}^{2}&+ x_{42}^{2}&+ 320 x_{10}&-16 x_{13} \\
-16 x_{14}&-16 x_{32}&-16 x_{33}&+ x_{1}^{4}&+ x_{2}^{4} \\
+ x_{3}^{2}&+ x_{4}^{2}&+ x_{5}^{2}&+ x_{7}^{2}&+ x_{9}^{2} \\
-2 x_{15}&-2 x_{34}&+ 704 x_{6}&+ 64 x_{1}&+ 2 x_{12} \\
+ 2 x_{16}&+ 2 x_{17}&+ 2 x_{19}&+ 2 x_{20}&+ 2 x_{21} \\
+ 2 x_{22}&+ 2 x_{23}&+ 2 x_{35}&+ 2 x_{36}&+ 2 x_{38} \\
+ 6 x_{41}&+ 4 x_{42}&-2 x_{5}&+ 6 x_{2}&+ 270 \\
 = 0 .
\end{array}
\end{align*}
}

\begin{remark}
This equation has a rather astonishing property, which is that for all $n \in \Z^+$, the number of solutions $\vec{x} \in \mathbb{N}^{42}$ is the $n$-th prime, $p(n)$.
\end{remark}
The exponential Diophantine equation $\hat{F}(n,\vec{x})=0$, where $\vec{x} = (x_1,x_2,\ldots,x_{42})$, contains $498$ monomials in $42$ quantified variables. In this expansion, the variables from \cref{eq:F1} have been relabeled over the tuple $\vec{x}$. So,
\begin{align*}
     & \hat{F}(n,\vec{x}) = \hat{F}(n,x_1,x_2,\ldots,x_{42}) = 0 \\
     & \iff
     \hat{F}(n,x_1,x_2,\ldots,x_8,[34])
     = \hat{F}(n,x_1,x_2,\ldots,x_8,y_1,y_2,\ldots,y_{34})
     = 0 \\
     & \iff
     \hat{F}(n, a, f_1, f_2, f_3, f_4, m, b, d, [34])
     = \hat{F}(n, a, f_1, f_2, f_3, f_4, m, b, d, y_1,y_2,\ldots,y_{34})
     = 0 .
\end{align*}
The SageMath source code used to generate the LaTeX for this expansion is included in Appendix \cref{appendix:sagesourcecodes} and can be examined to see exactly how this relabeling was performed.

\section{Maple source code for arithmetic terms} \label{appendix:maple}
The following Maple source code can be used to verify the arithmetic terms in \cref{proof:ntermarithmetic}, \cref{proof:omega}, \cref{proof:primepi}:
{\fontsize{8}{9}\selectfont
\begin{verbatim}
NU2 := n -> padic:-ordp(n, 2);
HW := n -> add(convert(n, base, 2));
GCD_arith := (m, n) -> irem(floor((2^(m^2*n*(n + 1)) - 2^(m^2*n))*(2^(m^2*n^2) - 1)
/((2^(m^2*n) - 1)*(2^(m*n^2) - 1)*2^(m^2*n^2))), 2^(n*m));
NU2_arith := n -> floor(irem(GCD_arith(n, 2^n)^(n + 1), (2^(n + 1) - 1)^2)/(2^(n + 1) - 1));
HW_airth := n -> nu2(irem(floor((1 + 2^(2*n))^(2*n)/2^(2*n^2)), 2^(2*n)));
G[0] := (q, t) -> (q^(t + 1) - 1)/(q - 1);
G[1] := (q, t) -> q*(t*q^(t + 1) - (t + 1)*q^t + 1)/(q - 1)^2;
G[2] := (q, t) -> q*(t^2*q^(t + 2) - (2*t^2 + 2*t - 1)*q^(t + 1)
+ (t + 1)^2*q^t - q - 1)/(q - 1)^3;
G[4] := (q, t) -> q*(t^4*q^(t + 4) + (-4*t^4 - 12*t^3 - 6*t^2 + 12*t + 11)*q^(t + 1)
+ (6*t^4 + 12*t^3 - 6*t^2 - 12*t + 11)*q^(t + 2)
+ (-4*t^4 - 4*t^3 + 6*t^2 - 4*t + 1)*q^(t + 3)
+ (t + 1)^4*q^t - q^3 - 11*q^2 - 11*q - 1)/(q - 1)^5;
C := (e, k, t, w) -> (2^w - e + 1)*(2^(2*w*t^k) - 1)/(2^w + 1);
A := proc(a, U, B, V, k, t, u)
local i;
return -(2^u - 1)*a*mul(G[U[i]](B[i]^V[i]*2^(2*u*t^(i - 1)), t - 1), i = 1 .. k);
end proc;
k := 2;
t := n -> n + 1;
u := n -> n + 5;
M := n -> C(1, 2, t(n), u(n))
+ A(1, [4, 0], [k, k], [0, 0], k, t(n), u(n))
+ A(-2, [2, 0], [k, k], [0, 0], k, t(n), u(n))
+ A(-2*n, [2, 1], [k, k], [0, 0], k, t(n), u(n))
+ A(n^2, [0, 2], [k, k], [0, 0], k, t(n), u(n))
+ A(2*n, [0, 1], [k, k], [0, 0], k, t(n), u(n));
N := n -> HW(M(n))/u(n) - t(n)^2;
omega := n -> NU2(N(4*n)) - 1;
seq(n, n = 1 .. 16);
seq(omega(n), n = 1 .. 16);
\end{verbatim}
}

The Maple source code above has been derived from the source code provided in \cite{prunescu2024numbertheoreticfunctions}. For presentation purposes, newline characters \textbackslash n have been inserted inside of the various arrow functions and procedures, which may cause compilation issues. To ensure the code executes properly, these functions and procedures should be rewritten to be on a single line.

\section{SageMath source code for arithmetic terms} \label{appendix:sagesourcecodesarithmeticterms}
The following SageMath source code can be used to verify the arithmetic terms in \cref{proof:ntermarithmetic}, \cref{proof:omega}, \cref{proof:primepi}:

{\fontsize{8}{9}\selectfont
\begin{verbatim}
from sage.all import *
def print_values(values, title=''):
    if title != '': print(f'{title}:')
    for v in values: print(f'{v}', end=',')
    print('')
def nu2(a): return a.valuation(2)
def HW(a): return bin(a).count('1')
def G0(q, t): return (q**t - 1) / (q - 1)
def G1(q, t):
    t1 = t - 1
    return q * (t1 * q**t - t * q**t1 + 1) / (q - 1)**2
def G2(q, t):
    t1 = t - 1
    return q * (
        t1**2 * q**(t1 + 2)
        - (t1**2 * 2 + t1 * 2 - 1) * q**(t1 + 1)
        + (t1 + 1)**2 * q**t1 - q - 1
    ) / (q - 1)**3
def G4(q, t):
    t1 = t - 1
    return q * (
        t1**4 * q**(t1 + 4)
        + (t1**4 * (-4) - 12 * t1**3 - 6 * t1**2 + t1 * 12 + 11) * q**(t1 + 1)
        + (t1**4 * 6 + 12 * t1**3 - 6 * t1**2 - t1 * 12 + 11) * q**(t1 + 2)
        + (t1**4 * (-4) - 4 * t1**3 + 6 * t1**2 - t1 * 4 + 1) * q**(t1 + 3)
        + (t1 + 1)**4 * q**t1 - q**3 - q**2 * 11 - q * 11 - 1
    ) / (q - 1)**5
def C(a, k, t, u): return ((2**u - a + 1) * (2**(u * 2 * t**k) - 1)) / (2**u + 1)
def G(r, q, t):
    if r == 0: return G0(q, t)
    elif r == 1: return G1(q, t)
    elif r == 2: return G2(q, t)
    elif r == 4: return G4(q, t)
    else:
        g = QQ(0)
        qj = QQ(1)
        for j in range(t):
            g += qj * j**r
            qj *= q
        return g
def A(a, P, B, V, k, t, u):
    result = -(2**u - 1) * a
    for i in range(k):
        q = B[i]**V[i] * 2**(u * t**i * 2)
        g = G(P[i], q, t)
        result *= g
    return result
# Initialize variables.
n_values = [i for i in range(1, 19)]
k = 2
n,t,u = var('n,t,u')
B = [2, 2]
V = [0, 0]
def M(a):
    return (
        C(1, k, t, u)
        + A(1, [4, 0], B, V, k, t, u)
        - A(2, [2, 0], B, V, k, t, u)
        - A(a*2, [2, 1], B, V, k, t, u)
        + A(a**2, [0, 2], B, V, k, t, u)
        + A(a*2, [0, 1], B, V, k, t, u))
def T(a): return a+1
def U(a): return a+5
# Construct the arithmetic term for M(n).
M_term = M(n)
def N(a):
    # Substitute values into the arithmetic term for M(n).
    p1 = M_term.subs(t=T(a)).subs(u=U(a)).subs(n=a)
    v1 = Integer(p1)
    v2 = HW(v1)/U(a)-T(a)^k
    return v2
def Omega(a): return nu2(N(a*4))-1
def PrimePi(a): return Omega(factorial(a))
# Display results.
print_values(n_values, 'N(n)')
print_values([N(a) for a in n_values])
print_values(n_values, 'Omega(n)')
print_values([Omega(a) for a in n_values])
n_values = [i for i in range(1, 5)]
print_values(n_values, 'Pi(n)')
print_values([PrimePi(a) for a in n_values])
\end{verbatim}
}

\section[SageMath source code for the n-th prime]{SageMath source code for the $n$-th prime equation} \label{appendix:sagesourcecodes}
The following code is used to generate the monomial expansion of the exponential Diophantine equation $\hat{F}(n,\vec{x})=0$ in Appendix \cref{appendix:F1monomials} and its corresponding arithmetic term.

{\fontsize{8}{9}\selectfont
\begin{verbatim}
from sage.all import *
enable_y_relablings = True          # Enables the y variable relabelings.
enable_32_variable_version = False  # Enables the 32 variable version equation F=0.
print_monomial_expansion = True     # Enables printing of LaTeX for \hat{F}=0.
print_details = True                # Enables printing of details about \hat{F}=0.
print_relabelings = True            # Enables printing of the relabelings for quantified variables.
enable_expand_cols = True           # Enables page breaks for generated LaTeX.
expand_cols_amount = 2      # Number of columns the array size for equations is expanded by
                            # following the first page break.
k = 42                      # Number of variables to initialize.
n = var('n')                # Placeholder variable for n.
x = var(['x{}'.format(i) for i in range(0, k+1)]) # This is our \vec{x} = (x_1,...,x_k) .
y = var(['y{}'.format(i) for i in range(0, k+1)]) # This is our \vec{y} = (y_1,...,y_k) .
relabelings = {}
def evaluate_polynomial(p, s, c):
    # Evaluates the polynomial p by substituting: s = c,
    # where s is the variable to be replaced.
    subs_dict = { s : c for i in range(1)}
    return p.subs(subs_dict)
def E_M4_9(x1,x2,x3,x4,i,x5):
    t = 4*x1+1              # t (non-quantified)
    u = 4*x1+5              # u (non-quantified)
    t1 = t-1                # t_1 (non-quantified)
    q11 = y[i+1]            # q_{1,1}
    q21 = y[i+2]            # q_{2,1}
    g01 = y[i+3]            # g_{0,1}
    g21 = y[i+4]            # g_{2,1}
    g41 = y[i+5]            # g_{4,1}
    g02 = y[i+6]            # g_{0,2}
    g12 = y[i+7]            # g_{1,2}
    g22 = y[i+8]            # g_{2,2}
    C = y[i+9]             # C
    h = (1-2^u)
    # In the polynomial expressions below, powers of q_1 and q_2 have been
    # explicitly replaced by their corresponding powers of 2 to ensure proper
    # usage of the relation: f_2 = f_1^2.
    q1 = 2^(2*u)            # q_1 (non-quantified)
    q2 = 2^(32*x2+48*x1+10) # q_2 (non-quantified)
    q1t = q2
    q2t = 2^(128*x3+352*x2+280*x1+50)
    q1t1 = 2^(32*x2+8*x1)
    q2t1 = 2^(128*x3+192*x2+40*x1)
    g411 = (-4)*t1^4 - 12*t1^3 - 6*t1^2 + 12*t1 + 11
    g412 = (6*t1^4 + 12*t1^3 - 6*t1^2 - 12*t1 + 11)
    g413 = (-4)*t1^4 - 4*t1^3 + 6*t1^2 - 4*t1 + 1
    g414 = (x4^4 * q1t1 + t1^4*q1t1 - q1^3 - 11*q1^2 - 11*q1 + 3)
    polys = [
        x5 - (C + h*g41*g02 - 2*h*g21*g02 - 2*(4*x1)*h*g21*g12 + (4*x1)^2*h*g01*g22
        + 2*(4*x1)*h*g01*g21),
        C*(2^u+1) - 2^u * (q2t - 1),
        q11 - q1 + 1,
        q21 - q2 + 1,
        g01*q11 - q21,
        g02*q21 - (q2t - 1),
        g12*q21^2 - q2*( t1 * q2t - x4 * q2t1 ),
        g21*q11^3 - q1*( t1^2 * q1t1*q1^2 - (2*t1^2+2*t1-1) * q2 + x4^2 * q1t1 - q1 - 1 ),
        g22*q21^3 - q2*( t1^2 * q2t1*q2^2 - (2*t1^2+2*t1-1) * q2t + x4^2 * q2t1 - q2 - 1 ),
        g41*q11^5 - q1*(q1t*g411 + q1t1*q1^2*g412 + q1t1*q1^3*g413 + g414 )
    ]
    relabelings.update({ 'q_{1,1}': q11, 'q_{2,1}': q21 })
    relabelings.update({ 'g_{0,1}': g01, 'g_{2,1}': g21, 'g_{4,1}': g41 })
    relabelings.update({ 'g_{0,2}': g02, 'g_{1,2}': g12, 'g_{2,2}': g22 })
    relabelings.update({ 'C': C })
    return build_poly(polys)
def E_divides_1(x1, i, x2):
    return (x1 - x2*y[i+1])^2
def E_notdivides_3(x1, i, x2):
    polys = [
        x1 - x2*y[i+1] - y[i+2] - 1,
        y[i+2] + y[i+3] + 2 - x2
    ]
    return build_poly(polys)
def E_nu_4(x1,i,x2):
    p1 = E_divides_1(x1, i, 2^x2)
    p1 += E_notdivides_3(x1, i+1, 2^(x2+1))
    return p1
def E_div_2(x1,x2,i,x3):
    polys = [
        x1-x2*x3-y[i+1],
        y[i+1]+y[i+2]+1-x3
    ]
    return build_poly(polys)
def E_mod_2(x1,x2,i,x3):
    polys = [
        x1-x2*y[i+1]-x3,
        x3+y[i+2]+1-x2
    ]
    return build_poly(polys)
def E_bin_7(x1,x2,i,x3):
    p1 = 0
    y1 = y[i+1]
    y2 = y[i+2]
    y3 = y[i+3]
    p1 += (y1-(2*x1^3+8*x1^2+2*x1*x2+12*x1+4*x2+8))^2
    p1 += (y2-(2*x1^2+8*x1+8))^2
    p1 += E_div_2(2^(y1), 2^(y2)-2^(2*x1+4)-1, i+3, y3)
    p1 += E_mod_2(y3, 2^(2*x1+4), i+5, x3)
    return p1
def E_fact_13(x1,i,x2):
    p1 = 0
    y1 = y[i+1]
    y2 = y[i+2]
    y3 = y[i+3]
    y4 = y[i+4]
    p1 += (y1 - x1^2)^2
    p1 += (y2 - 2^(3*x1))^2
    p1 += (y3 - x1*y1)^2
    p1 += E_bin_7(y2, x1, i+4, y4)
    p1 += E_div_2(2^(3*y3), y4, i+11, x2)
    return p1
def E_HW_12(x1,i,x2):
    p1 = 0
    y1 = y[i+1]
    p1 += E_bin_7(2*x1, x1, i+1, y1)
    p1 += E_nu_4(y1, i+8, x2)
    return p1
def build_poly(polys):
    p1 = 0
    for p in polys: p1 += p^2
    return p1
def is_constant_monomial(p_operand):
    e1 = evaluate_polynomial(p_operand, n, 1)
    if e1.is_constant():
        return True
    if e1.operator() == operator.pow:
        base, exponent = e1.operands()
        if base.is_constant() and exponent.is_constant():
            return True
    return False
def G0(q, t): return (q**t - 1) / (q - 1)
def G1(q, t):
    t1 = t - 1
    return q * (t1 * q**t - t * q**t1 + 1) / (q - 1)**2
def G2(q, t):
    t1 = t - 1
    return q * (
        t1**2 * q**(t1 + 2)
        - (t1**2 * 2 + t1 * 2 - 1) * q**(t1 + 1)
        + (t1 + 1)**2 * q**t1 - q - 1
    ) / (q - 1)**3
def G4(q, t):
    t1 = t - 1
    return q * (
        t1**4 * q**(t1 + 4)
        + (t1**4 * (-4) - 12 * t1**3 - 6 * t1**2 + t1 * 12 + 11) * q**(t1 + 1)
        + (t1**4 * 6 + 12 * t1**3 - 6 * t1**2 - t1 * 12 + 11) * q**(t1 + 2)
        + (t1**4 * (-4) - 4 * t1**3 + 6 * t1**2 - t1 * 4 + 1) * q**(t1 + 3)
        + (t1 + 1)**4 * q**t1 - q**3 - q**2 * 11 - q * 11 - 1
    ) / (q - 1)**5
def C(a, k, t, u): return ((2**u - a + 1) * (2**(u * 2 * t**k) - 1)) / (2**u + 1)
def G(r, q, t):
    if r == 0: return G0(q, t)
    elif r == 1: return G1(q, t)
    elif r == 2: return G2(q, t)
    elif r == 4: return G4(q, t)
    else:
        g = QQ(0)
        qj = QQ(1)
        for j in range(t):
            g += qj * j**r
            qj *= q
        return g
def A(a, P, B, V, k, t, u):
    result = -(2**u - 1) * a
    for i in range(k):
        q = B[i]**V[i] * 2**(u * t**i * 2)
        g = G(P[i], q, t)
        result *= g
    return result
k = 2
t,u = var('t,u')
B = [2, 2]
V = [0, 0]
def M(a):
    return (
        C(1, k, t, u)
        + A(1, [4, 0], B, V, k, t, u)
        - A(2, [2, 0], B, V, k, t, u)
        - A(a*2, [2, 1], B, V, k, t, u)
        + A(a**2, [0, 2], B, V, k, t, u)
        + A(a*2, [0, 1], B, V, k, t, u))
if enable_y_relablings: y = x    # Relabel \vec{y} variables over \vec{x}.
i = 0
F = 0
a = x[i+1]
f1 = x[i+2]
f2 = x[i+3]
f3 = x[i+4]
f4 = 0
m3 = 0
if enable_32_variable_version:
    m = x[i+5]
    b = x[i+6]
    d = x[i+7]
    i += 7
    F += E_fact_13(a, i, f1)      # Add the sum of squares for the relation: f_1 = a! .
    i += 13
    F += (f2 - f1^2)^2            # Add the square for the relation: f_2 = f_1^2 .
    F += (f3 - f1*f2)^2           # Add the square for the relation: f_3 = f_1^3 .
    f4 = 4*f1+1
    m0 = M(4*f1).subs(t=4*f1+1).subs(u=4*f1+5)
    # NOTE: We are not replacing f_1^2=f_2 and f_1^3=f_3 in the exponents here because
    # we only care about the monomial counts, which are the same regardless of replacement.
    # To obtain the 32 variable monomial expansion which works with the hypercube method,
    # one must update the code to replace these exponents.
    m1 = m0.simplify_rational()
    D = m1.denominator()
    L = m1.numerator()
    m3 = (m*D-L)^2
    F += m3                      # Add the square for the relation: m = M(4(a!)).
else:
    f4 = x[i+5]
    m = x[i+6]
    b = x[i+7]
    d = x[i+8]
    i += 8
    F += E_fact_13(a, i, f1)      # Add the sum of squares for the relation: f_1 = a! .
    i += 13
    F += (f2 - f1^2)^2            # Add the square for the relation: f_2 = f_1^2 .
    F += (f3 - f1*f2)^2           # Add the square for the relation: f_3 = f_1^3 .
    F += (f4 - (4*f1+1))^2              # Add the square for the relation: f_4 = 4 f_1 + 1 .
    F += E_M4_9(f1, f2, f3, f4, i, m)   # Add the sum of squares for the relation: m = M(4(a!)) .
    i += 9
F += E_HW_12(m, i, b)   # Add the sum of squares for the relation: b = HW(m) .
i += 12
F += (b + (f1+5)*(-(f4^2)+d-2^n))^2  # Add the square for the relation: 2^{\pi(a)+1} <= 2^n .
relabelings.update({ 'a': a, 'f_1': f1, 'f_2': f2, 'f_3': f3, 'f_4': f4, 'b': b, 'd': d })
F = F.expand()
monomials = F.operands()
constants = [po for po in monomials if is_constant_monomial(po)]
nonconstants = [po for po in monomials if is_constant_monomial(po) == False]
p0_constants = sum(constants)
p0_nonconstants = sum(nonconstants)
p0_constants_text = latex(p0_constants)
if print_details:
    print(f'Details: Num. Variables = {len(F.variables())-1}, Num. Monomials = {len(monomials)}')
    if enable_32_variable_version:
        print(f'M Monomials: {len(m3.expand().operands())}')
    print(f'Num. ATerms: {len(nonconstants)}.')
    print(f'Constants:\n{p0_constants_text}')
if print_relabelings:
    relabelings_text = ''
    for key,value in relabelings.items():
        relabelings_text += f'{key} = {latex(value)}\n'
    print(f'Quantified variable relabelings:\n{relabelings_text}')
    print(f'NOTE: The remaining relabelings can be found by examining the functions')
    print(f'E_fact_13(.), E_M_9(.), and E_HW_12(.) in the order that they are invoked.')
def format_operand_monomial(po):
    po_text = latex(factor(po)).replace('\\, ','')
    return po_text
def get_operands_text(p_operands, format_func, arraylen=3, j=1, page_lines=60, num_lines=0):
    text = ''
    ltext = 'l'*(arraylen)
    num_lines = 1
    array_cols_expanded = False
    if j != 1:
        num_lines += j // arraylen
    op_strings = []
    for po in p_operands:
        po_text = format_func(po)
        op_strings.append(po_text)
    sorted_op_strings = sorted(op_strings, key=lambda x: len(x), reverse=True)
    is_new_line = False
    is_page_break = False
    num_page = 1
    for po_text in sorted_op_strings:
        if j != 1 or is_page_break: text += '+ ' if is_new_line else '&+ '
        is_new_line = False
        is_page_break = False
        if j != 1 and j % arraylen == 0:
            is_new_line = True
            num_lines += 1
            if num_lines % page_lines == 0:
                is_page_break = True
                if enable_expand_cols:
                    arraylen += expand_cols_amount*num_page
                    ltext = 'l'*(arraylen)
                    array_cols_expanded = True
                num_page += 1
        text += po_text
        j += 1
        if j <= len(sorted_op_strings):
            if is_new_line:
                if is_page_break:
                    text += '+ \n'
                    text += f'\\end{{array}}\n'
                    text += f'\\end{{align*}}\n'
                    text += f'\\begin{{align*}}\n'
                    text += f'\\begin{{array}}{{{ltext}}}\n'
                    j = 1
                else:
                    text += ' \\\\\n'
    return text
def get_latex(p_ops, format_func, prefix='', postfix='', arraylen=3, j=1, page_lines=70, num_lines=0):
    ltext = 'l'*(arraylen)
    text = ''
    text += f'\\begin{{align*}}\n'
    text += f'\\begin{{array}}{{{ltext}}}\n'
    if prefix != '':
        text += prefix
    operands_text = get_operands_text(p_ops, format_func, arraylen, j, page_lines)
    text += operands_text
    if postfix != '': text += postfix
    text += f'\n\\end{{array}}\n'
    text += f'\\end{{align*}}\n'
    text = text.replace('+ -', '-').replace('\\, ','')
    return text
if print_monomial_expansion:
    prefix_text = f'F(n,\\vec{{x}}) = ' if enable_32_variable_version else f'\\hat{{F}}(n,\\vec{{x}}) = '
    postfix_text = ' = 0 .'
    expand_cols_amount = 1
    num_cols = 2
    ptext = get_latex(monomials, format_operand_monomial, prefix_text, postfix_text, num_cols, 1, 55)
    print(f'Monomial Expansion:')
    print(ptext)
\end{verbatim}
}
\end{appendices}


\end{document}